\documentclass{amsart} \usepackage[mathscr]{eucal}
\usepackage{amssymb, amsmath,array, amscd}
\usepackage{enumerate}

\usepackage{hyperref}

\usepackage{xy}

\xyoption{all}

\def\aa{\ensuremath{\alpha}}

\def\wed#1{\ensuremath{\stackrel{#1}{\wedge}\!}}
\def\id#1{\ensuremath{ \langle {#1} \rangle}}
\def\p#1{\ensuremath{\mathbb{P}^{#1}}}
\def\ps#1{\ensuremath{\mathbb P\left( {#1}
\right)}}
\def\pd#1{\ensuremath{\check{\mathbb{P}}^{#1}}}
\def\C{\ensuremath{\mathbb C}}

\def\ra{\ensuremath{\rightarrow}}
\def\rar{\ensuremath{\dashrightarrow}}
\def\lar{\ensuremath{\longrightarrow}}
\def\ba#1{\begin{array}{#1}}\def\ea{\end{array}}
\def\be{\begin{equation}}\def\ee{\end{equation}}
\def\inj{\ensuremath{\xymatrix{\ar@{>->}[r]&}}}
\def\surj{\ensuremath{\xymatrix{\ar@{->>}[r]&}}}
\def\emb{\raise-1pt\hbox{\ensuremath{\xymatrix{\ar@{^(->}[r]&}}}}
\def\cl#1{\ensuremath{\mathcal{#1}}}
\def\bbb#1{\ensuremath{\mathbb{#1}}}
\def\O{\cl O}
\def\G{\ov{\mathbf S}\ensuremath{_{d_1}}}
\def\us{\ensuremath{^\star}}
\def\ov#1{\ensuremath{\overline{#1}}}
\def\nn{\ensuremath{\ne}} \def\ie{{\em i.e.\em,
}} \def\eg{\em e.g., \em}
\def\rank{\operatorname{rank}}
\def\Sing#1{\operatorname{sing}(#1)}
\def\cod{\operatorname{codim}}

\def\Hom{{\operatorname {Hom}}}
\def\ds#1{\displaystyle{#1}}
\def\bsm{\left(\begin{smallmatrix}}\def\esm{\end{smallmatrix}\right)}
\def\na#1{\noalign{\vskip#1pt}} \def\vez
{\ensuremath{\times}} \def\u#1{\underline{#1}}
\def\wt#1{\ensuremath{\widetilde{#1}}}
\def\rf#1{\hbox{\rm(\ref{#1})}} \def\eg{\em
e.g., \em} \def\inv{\ensuremath{{}^{\!-1}}}
\def\ls{\ensuremath{_\star}}
\def\bi{\begin{itemize}}\def\ei{\end{itemize}}
\def\Sym{\operatorname{Sym}} \def\n{\noindent}

\def\sing{{\sc Singular}\,\cite{sing}}
\def\schub{{\sc schubert}\,\cite{schub}}

\def\st{\ensuremath{\,|\,}}
\def\h{\ensuremath{\mathbf h}}
\def\q{\ensuremath{\mathbf q}}
\def\s#1{\ensuremath{\mathbf S_{#1}}}
\def\EE{\ensuremath{\mathbf E}}
\def\E{\ensuremath{\mathbf e}}
\def\B{\ensuremath{\mathbf B}}
\def\V{\ensuremath{\mathbf V}}
\def\x{\ensuremath{\mathbf X}}

 \def\aa{\ensuremath{\alpha}}

\def\wed#1{\ensuremath{\stackrel{#1}{\wedge}\!}}
\def\id#1{\ensuremath{ \langle {#1} \rangle}}
\def\p#1{\ensuremath{\mathbb{P}^{#1}}}
\def\ps#1{\ensuremath{\mathbb P\left( {#1} \right)}}
\def\pd#1{\ensuremath{\check{\mathbb{P}}^{#1}}}
\def\C{\ensuremath{\mathbb C}}  \def\ra{\ensuremath{\rightarrow}}
\def\rar{\ensuremath{\dashrightarrow}}
\def\lar{\ensuremath{\longrightarrow}}
\def\ba#1{\begin{array}{#1}}\def\ea{\end{array}}
\def\be{\begin{equation}}\def\ee{\end{equation}}
\def\inj{\ensuremath{\xymatrix{\ar@{>->}[r]&}}}
\def\surj{\ensuremath{\xymatrix{\ar@{->>}[r]&}}}
\def\emb{\raise-1pt\hbox{\ensuremath{\xymatrix{\ar@{^(->}[r]&}}}}
\def\cl#1{\ensuremath{\mathcal{#1}}}
\def\bbb#1{\ensuremath{\mathbb{#1}}}
\def\O{\cl O}
\def\G{\ov{\mathbf S}\ensuremath{_{d_1}}}
\def\us{\ensuremath{^\star}}
\def\ov#1{\ensuremath{\overline{#1}}}
\def\nn{\ensuremath{\ne}} \def\ie{{\em i.e.\em, }}
\def\eg{\em e.g., \em} \def\rank{\operatorname{rank}}
\def\Sing#1{\operatorname{sing}(#1)}
\def\cod{\operatorname{codim}}  \def\Hom{{\operatorname {Hom}}}
\def\ds#1{\displaystyle{#1}}
\def\bsm{\left(\begin{smallmatrix}}\def\esm{\end{smallmatrix}\right)}
\def\na#1{\noalign{\vskip#1pt}} \def\vez
{\ensuremath{\times}} \def\u#1{\underline{#1}}
\def\wt#1{\ensuremath{\widetilde{#1}}}
\def\rf#1{\hbox{\rm(\ref{#1})}} \def\eg{\em e.g., \em}
\def\inv{\ensuremath{{}^{\!-1}}}
\def\ls{\ensuremath{_\star}}
\def\bi{\begin{itemize}}\def\ei{\end{itemize}}
\def\Sym{\operatorname{Sym}} \def\n{\noindent}

 \def\sing{{\sc
Singular}\,\cite{sing}} \def\schub{{\sc
schubert}\,\cite{schub}}

\def\st{\ensuremath{\,|\,}}
\def\h{\ensuremath{\mathbf h}}\def\q{\ensuremath{\mathbf q}}
\def\s#1{\ensuremath{\mathbf S_{#1}}}
\def\v{\mathrm v}
\def\ve{\ensuremath{^\vee}}
\newtheorem{thm}{Theorem}[section]
 \newtheorem{THM}{Theorem}

\newtheorem{prop}{Proposition}[section]

\newtheorem{lemma}{Lemma}[section]

\newtheorem{remark}{Remark}[section]
\newtheorem{cor}{Corollary}[section]

\begin{document}
\title[Stability of foliations]{Stability of foliations induced by  rational maps}
\author{F. Cukierman} \author{J. V. Pereira}
\author{I. Vainsencher}

\thanks{The authors were partially supported by CAPES-SPU.
The second author is supported by Instituto Unibanco and CNPQ. }



\begin{abstract}
We show that the singular holomorphic foliations
induced by dominant quasi-homogeneous rational maps
fill out irreducible components of the  space $\mathscr F_q(r, d)$
of singular  foliations of codimension $q$ and degree $d$
on the complex projective space  $\mathbb P^r$, when $1\le  q \le r-2$.
We study the geometry of these irreducible components.
In particular we prove that they are all rational varieties
and we compute their projective degrees in several cases.
\end{abstract}

\thanks{} \keywords{} \subjclass{}

\maketitle

\section{Introduction}
\label{section Introduction}

\subsection{The space of codimension one holomorphic foliations on $\mathbb P^r$}

\medskip

Let us consider a differential 1-form in $\mathbb C^{r+1}$
$$ \omega = \sum_{i=0}^r a_i dx_i
$$ where the $a_i$ are homogeneous polynomials of degree
$d+1$ in variables $x_0, \dots, x_r$, with complex
coefficients. Assume that $ \sum_{i=0}^r a_i x_i = 0$, so
that $\omega$ descends to the complex projective space $\p
r$ and defines a global section of the twisted sheaf of
1-forms $\Omega^1_{\mathbb P^r}(d+2)$.

The space of codimension one foliations of degree $d$ on
$\mathbb P^r$ is the algebraic subset of $\ps {\mathrm
H^0(\p r, \Omega^1_{\p r}(d+2))}$ consisting of the 1-forms
$\omega$ that satisfy the Frobenius integrability condition
and has zero set of codimension at least two, i.e.,
$$ \mathscr F(r, d) = \left\{ \omega \in \ps {\mathrm H^0(\p
r, \Omega^1_{\p r}(d+2))} \, | \, \omega \wedge d \omega =0
\text{ and } \cod \Sing{\omega} \ge 2 \right\}.
$$

For the study of the irreducible components of $\mathscr F(r,d)$
we refer to e. g.  \cite{CL}   and  \cite{Jouanolou}.

\subsection{Stability of quasi-homogeneous  pencils}

\medskip

One of the first results on the subject is due to G\'{o}mez-Mont and
Lins Neto \cite{GMLN} who proved that there are irreducible
components $\mathscr R(r,d,d) \subset \mathscr F(r,2d-2)$, $r\ge 3$,
whose generic element is a foliation tangent to a Lefschetz pencil
of degree $d$ hypersurfaces. Their proof explores the topology of
the underlying real foliation and relies on the stability of the
Kupka components of the singular set and on Reeb's Leaf Stability
Theorem. Using similar methods they recognized for $r\ge 4$ other
irreducible components $\mathscr R(r,d_0,d_1) \subset \mathscr F(r,
d_0 + d_1 -2)$ with generic member tangent to a quasi-homogeneous
pencil $\id{\lambda F^{p_0} - \mu G^{p_1}}$ with $p_0$ and $p_1$
relatively prime natural numbers satisfying $p_0 d_0 =p_1 d_1,
\,d_i=\deg F_i$. Later Calvo-Andrade  \cite{Omegar} extended
G\'{o}mez-Mont-Lins Neto result about quasi-homogeneous pencils to
dimension three. His proof has an extra dynamical ingredient --the
stability of leaves carrying non-trivial holonomy.

In fact in both of the above mentioned papers the authors do
not restrict to $\mathbb P^r$ and prove their results for
foliations on an arbitrary projective manifold $M$ with
$\dim M\ge 3$ and $\mathrm H^1(M,\mathbb C)=0$. An
alternative proof of the above results based on extension
techniques of transversely euclidean structures has been
carried out by Sc\'{a}rdua in \cite{SC}.

\subsection{Infinitesimal stability of quasi-homogeneous  pencils}

\medskip

Although full of geometric insights the above mentioned
works do not seem to shed any light on the scheme structure
or the geometry of $\mathscr R(r,d_0,d_1)$. The present
article stems from an attempt to understand these problems.

Using infinitesimal techniques, as in \cite{CP}, we describe
the Zariski tangent space of $\mathscr R(r,d_0,d_1)$ at a
generic point and arrive at a proof that $\mathscr R(r, d_0,
d_1)$ --with the natural scheme structure given by the
Frobenius integrability condition-- is {\it generically
reduced}. More precisely if $\mathscr R(r, d_0, d_1)$
denotes the closure of the image of the rational map
\begin{eqnarray*}
\rho: \ps{ \mathrm H^0(\p r, \mathscr O_{\p r}(d_0))} \times
\ps{\mathrm H^0(\p r, \mathscr O_{\p {r}}(d_1))} &
\mbox{\Large$\dashrightarrow$}& \ps{\mathrm H^0(\p
r,\Omega^1(d_0 + d_1))} \\ (F_0, F_1)
&\mbox{\Large$\mapsto$} & d_0 F_0 dF_1 - d_1 F_1 dF_0.
\end{eqnarray*}
then our first result reads as follows.

\begin{THM}\label{T:1}
If $r\ge 3$ then $\mathscr R(r, d_0, d_1)$ is an irreducible
and \underline{generically reduced} component of $\mathscr
F(r, d_0 + d_1 -2 )$.
\end{THM}

As explained above the only novelty in Theorem \ref{T:1},
besides the method of its proof, is what concerns the scheme structure
over a generic point. For a more precise statement see
Theorem \ref{T:reduzida} in \S \ref{S:cod1}.

The main content of this article is the generalization of
Theorem \ref {T:1} to foliations of higher codimension.

\subsection{Foliations on $\mathbb P^r$ of higher codimension}

\medskip

Let $\omega$ be a homogeneous $q$-form on $\mathbb C^{ r+1}$
with coefficients of degree $d+1$ that is annihilated by
Euler's vector field. As before $\omega$ can be interpreted
as a section of the sheaf of twisted differential $q$-forms
$\Omega^q_{\p r}(d+q+1)$.

We recall from \cite{Airton} (see also \cite{CP}) that
$\omega$ defines a degree $d$ holomorphic foliation of
codimension $q$ on $\mathbb P^r$ if it satisfies both
Pl\"{u}cker's decomposability condition
\begin{equation}\label{E:pluckeri}
 ( {i}_v\omega) \wedge \omega = 0 \quad \text{ for every }v
 \in \bigwedge^{q-1}\mathbb C^{r+1} ,
\end{equation}
and the integrability condition
\begin{equation}\label{E:pluckerii}
 ( {i}_v\omega) \wedge d\omega = 0 \quad \text{ for every }v
 \in \bigwedge^{q-1}\mathbb C^{r+1} .
\end{equation}

It is therefore natural to set $\mathscr F_q(r,d)$, the
space of codimension $q$ holomorphic foliations of degree
$d$ on $\mathbb P^r$, as
\[
  \left\lbrace \omega \in \ps {\mathrm H^0(\p r,
\Omega^q_{\p r}(d+q+1))} \, \big| \, \omega \text{ satisfies
} (\ref{E:pluckeri}),(\ref{E:pluckerii}) \text{ and }
\mathrm{codim} \ \mathrm{sing}(\omega) \ge 2 \right\rbrace .
\]

\subsection{Infinitesimal stability of quasi-homogeneous rational maps}

\medskip

If one interprets the elements of $\mathscr R(r,d_0,d_1)$ as
foliations tangent to the fibers of rational maps
\begin{eqnarray*}
  \mathbb P^r &\dashrightarrow& \mathbb P^1 \\ x &\mapsto& (
 F^{p_0} : G^{p_1} )
\end{eqnarray*}
then a possible counterpart in the higher codimension case are
the foliations tangent to dominant rational maps
$\mathbb P^r \dashrightarrow \mathbb P^q$.

When $q= r-1$ there is no hope to establish a stability
result even for a generic rational map. Indeed, under this
constraint both Pl\"{u}cker's condition and the
integrability condition are vacuous.  Thus $\mathscr
F_{r-1}(r,d)$ can be identified with an open subset of
$\ps{\mathrm H^0(\p r, \Omega^{r-1}_{\p r}(d+r))} =
\ps{\mathrm H^0(\p r, T \p r(d-1))}$. It is well known that
for $d\ge 2$ a generic element of this space has no
algebraic leaves, see for instance \cite{Collier}.

For $1 \le q \le r-2$ fix integers $d_0, \dots,d_q$ and
consider homogeneous polynomials $F_i$ of degree $d_i$ for
$i = 0, \dots, q$.  Assume that the $q$-form
\begin{equation}\label{E:dq}
\omega = i_R(dF_0 \wedge \dots \wedge dF_q),
\end{equation}
is non-zero.
It is easy to check that $\omega$
satisfies both (\ref{E:pluckeri}) and
(\ref{E:pluckerii}) since
$i_v \omega= \sum  a_{ij} i_R (dF_i \wedge dF_j)$,
where the $a_{ij}$ are homogeneous polynomials.
Moreover, it defines a foliation
tangent to the fibers of the map
\begin{eqnarray*}
  \mathbb P^r &\dashrightarrow& \ \mathbb P^q \\
x\ &\mapsto& (
  F_0^{e_0}: \ldots : F_q^{e_q} )
\end{eqnarray*}
with $e_i = \mathrm{lcm}(d_0,\ldots, d_q) /d_i$.
We set
$$
d=\sum d_i - q -1
$$
and denote by
$$
\mathscr R(r,d_0,\ldots, d_q) \subset \mathscr F_q(r,d)
$$
the closure of the set of foliations that can be written
in the form (\ref{E:dq}). It is the closure of the image of
the rational map
$$\ba{ccc}
\rho:\prod_i \ps{\mathrm H^0(\O_{\p r}(d_i))}&\dashrightarrow&
\ps{\mathrm H^0(\p r,\Omega^1(d+q+1))}\\
(F_i)&\mapsto&i_R(dF_0 \wedge \dots \wedge dF_q).
\ea
$$Notice that for $q=1$ we recover
the definition of $\mathscr R(r,d_0,d_1)$.

\begin{THM} \label{T:2}
If $r\ge 4$ and $1 \le q \le r-2$ then $\mathscr R(r, d_0,
\ldots, d_q)$ is an irreducible and generically reduced
component of $\mathscr F_q \left(r, \sum d_i -q -1 \right)$ .
\end{THM}

As far as we know there is no information in the literature
concerning the geometry of the irreducible components of
$\mathscr F_q(r,d)$ so far.

\subsection{Geometry of the rational components}

\medskip

 In Section \ref{S:par} we initiate this study through an
investigation of the parameterization $\rho$.
Besides computing the  dimension of $\mathscr R(r, d_0, \ldots, d_q)$,
we prove the following.

\begin{THM}\label{T:rational}
The {irreducible components}
$\mathscr R(r, d_0, \ldots, d_q)$ are {rational varieties}.
\end{THM}

By its definition,  $\mathscr R(r, d_0, \ldots, d_q)$ is  unirational.
The proof of rationality relies on the
construction of a variety $X$ that sits as an open set in the total space
of a tower of Grasmmann bundles, together
with a  birational morphism
$p: X \to \mathscr R(r,d_0,\ldots,d_q)$.

In general we do not know how to naturally compactify $X$ to
a projective variety where $p$ extends to a
morphism. Albeit, in a number of cases we are able to do
that and obtain, with the aid of Schubert Calculus, formulas
for the degree of the projective subvarities
$$ \mathscr R(r, d_0, \dots, d_q) \subset
\ps{\mathrm{H}^0(\p r,\Omega^q(d+q+1))} .
$$ For example the first few values for the degree of
$\mathscr R(r, 2, 2,2) $ are listed below. \medskip

\centerline{
\begin{tabular}{c}
\begin{tabular}{||l|l||}
\hline $r$&{Degree}\\\hline 3&{\bf1324220}\\\hline
4&{\bf2860923458080}\\\hline
5&{\bf243661972980477736263}\\\hline 6&
{\bf728440733705107831789517245858}\\\hline 7&
{\bf704613096513585123585398408696231899176183} \\\na{-2}
\hline\hline
\end{tabular}
\end{tabular}
}

\medskip

Several other cases are treated in Section \ref {S:section Degree}.

\tableofcontents

\section{Infinitesimal stability of quasi-homogeneous
pencils}\label{S:cod1}

In this first section we present our proof of Theorem
\ref{T:1}. All the arguments will be reworked later in
greater generality. We felt the exposition of this
particular case of Theorem \ref{T:2} would improve the
clarity of the paper.

\medskip

 For simplicity, let us denote by
\begin {equation} \label{polynomials}
\s e = \mathrm H^0(\p r, \mathscr O_{\p r}(e))
\end {equation}
the vector space of homogeneous polynomials of degree $e$ in
$r+1$ variables, and
$$ \mathscr F = \mathscr F(r, d)
$$ so that our rational map $\rho$ is
\begin {equation} \label{map}
\rho: \ps{ \s {d_0} } \times \ps{ \s {d_1} }
\dashrightarrow \mathscr F \subset \ps{ \mathrm H^0(\p r,
\Omega^1(d+2)) }.
\end {equation}

If $p_0$ and $p_1$ denote the unique coprime natural numbers
such that $ p_0 d_0 = p_1 d_1$ then
$$ \rho(F_0, F_1) = d_0 F_0 dF_1 - d_1 F_1 dF_0 = p_1 F_0
dF_1 - p_0 F_1 dF_0
$$ where the last equality of differential forms is up to
multiplicative constant.

We remark that
$$d \left(\frac{F_0^{p_0}}{F_1^{p_1}} \right) =
\frac{F_0^{p_0 - 1}}{F_1^{p_1 + 1}} \ (p_1 F_0 dF_1 - p_0
F_1 dF_0).
$$ Therefore, the closure of the leaves of the singular
foliation defined by the integrable 1-form $\rho(F_0, F_1)$
are irreducible components of the members of the pencil of
hypersurfaces of degree $p_0 d_0 = p_1 d_1$ generated by
$F_0^{p_0}$ and $F_1^{p_1}$.

\subsection{The Zariski tangent space of $\mathscr F$}

\medskip

 For a scheme $X$ and a point $x \in X$ we denote by $T_xX$
the Zariski tangent space of $X$ at $x$. If $\ps{V}$ is the
projective space associated to a $\mathbb C$-vector space
$V$ and denoting $\pi: V -\{0\} \to \ps{V}$ the canonical
projection, for each $v \in V$ we have a natural
identification
$$
T_{\pi(v)} \ps{V} = V /(v)
$$
where $(v)$ denotes de one-dimensional subspace generated
by $v$.  With slight abuse of notations, the Zariski tangent
space $T_{\omega} \mathscr F $ of $\mathscr F$ at a point
$\omega$ is represented by the forms $\eta \in \mathrm
H^0(\p r, \Omega^1(d+2)) / (\omega) $ such that
$$
(\omega + \epsilon \eta)\wedge (d\omega + \epsilon d
\eta) = 0 \mod \epsilon^2
$$ that is, such that
$$
\omega \wedge d\eta + \eta \wedge d \omega = 0 \, \, \,
\text{ or, equivalently }  \, \,\, \,\, d\omega
\wedge d\eta =0 ,
$$
where the equivalence is implied by the following variant
of Euler's formula for homogeneous polynomials.
\begin{lemma}\label{L:1}
If $\eta$ is a homogeneous $q$-form with degree $d$
coefficients then
$$
i_R d\eta + d( i_R \eta) = (q+d) \eta
$$
where $R$ is the radial or Euler vector field and $i_R$
denotes the interior product or contraction with $R$.
\end{lemma}
\begin{proof}
See \cite[Lemme 1.2, pp. 3]{Jouanolou}.
\end{proof}

Therefore to determine $T_\omega \mathscr F$ is equivalent
to solve $d \omega \wedge d \eta = 0$. Notice that in the
situation under scrutiny $d \omega = (d_0 + d_1) dF_0 \wedge
dF_1$. The first step towards the general $\eta$ satisfying
$d \omega \wedge d \eta = 0$ is given by Saito's generalization
of DeRham's division Lemma. In Lemma \ref{L:2} we state
variants of both DeRham's and Saito's Lemmas fine tuned up
for our purposes.

\begin{lemma}[\cite{Saito}]\label{L:2}
Let $F_0 , \ldots, F_q$ be homogeneous polynomial functions
on $\mathbb C^{r+1}$ and let $\Theta \in
\Omega^{q+1}(\mathbb C^{r+1})$ be the $(q+1)$-form given by
\[
\Theta = dF_0 \wedge \ldots \wedge dF_q \, .
\]
\begin{enumerate}
\item[{\rm (a)}] Suppose that $q< r$ and $\mathrm{codim } \,
\mathrm{ sing} (\Theta) \ge 2$. If $\eta \in
\Omega^1(\mathbb C^{r+1})$ is a homogeneous polynomial
$1-$form such that $\Theta \wedge \eta =0$ then there exist
homogeneous polynomials $a_0, \ldots, a_q$ such that
$$ \eta = \sum_{i=0}^q a_i  dF_i.
$$
\item[{\rm (b)}] Suppose that $q< r-1$ and $\mathrm{codim } \,
\mathrm{ sing} (\Theta) \ge 3$. If $\eta \in
\Omega^2(\mathbb C^{r+1})$ is a homogeneous polynomial
$2-$form such that $\Theta \wedge \eta =0$ then there exist
homogeneous polynomial $1$-forms $\alpha_0, \ldots,
\alpha_q$ such that
$$ \eta = \sum_{i=0}^q \alpha_i \wedge dF_i.
$$
\end{enumerate}
\end{lemma}

\begin{remark}\rm The hypothesis $q<r$ in (a) and $q < r-1$ in (b) are  not
really necessary.  For instance in item (b) the singular set
$\mathrm {sing }(\Theta)$ equals the locus where the $(q+1)
\times (r+1)$ Jacobian matrix $(\partial F_i / \partial
x_j)$ has rank $\le q$. Hence $\mathrm {sing }( \Theta)$ is
empty or has codimension at most $r+1-q$. When $q \ge r-1$
it follows that $\mathrm{codim } \, \mathrm{ sing} (\Theta)
\ge 3$ implies that $\Theta$ has no singularities. We
conclude that $F_0, \ldots, F_q$ are linearly independent
linear forms and the conclusion trivially holds true in this
case.
\end{remark}

In face of Lemma \ref{L:2} it is natural to define the
open subset
\begin{equation}\label{U} \mathcal U = \{ \omega \in \mathscr R(r, d_0,
d_1) \ | \ \mathrm{codim \ sing}(d\omega )\ge 3 \mathrm{\ and \ }
\mathrm{codim \ sing}(\omega)\ge 2 \}. \end{equation} The next
result will imply the infinitesimal stability of quasi-homogeneous
pencils corresponding to points of $\mathcal U$. It is a simple
particular case of Proposition \ref{P:derivativeq}. The iteration
argument in the proof is generalized in Lemma \ref{again}. We feel
it is worthwhile to write it here for the sake of clarity.

\begin{prop}\label{P:derivative} Let $(F_0, F_1) \in \ps{ \s {d_0}
} \times \ps{\s {d_1} }$ be such that $\rho(F_0, F_1)
= \omega \in \mathcal U$. Then the derivative
$$d \rho (F_0, F_1) : T_{(F_0,F_1)} ( \ps{ \s {d_0} }
\times \ps{\s {d_1} })  \to T_{\omega} \mathscr F$$ is
surjective. In other words, $\rho$ is a submersion over
$\mathcal U$.
\end{prop}
\begin{proof}
It is convenient to write
$$ \rho(F_0, F_1) = d_0 F_0 dF_1 - d_1 F_1 dF_0 = i_R (dF_0
\wedge dF_1) .
$$ Then, the derivative of $\rho$ at the point $(F_0,
 F_1)$
$$ d \rho (F_0, F_1): \s {d_0} / (F_0) \times \mathbf
S_{d_1} / (F_1) \rightarrow T_{\omega} \mathscr F
$$ is calculated as $$ d \rho (F_0, F_1) (F'_0, F'_1) = i_R
( dF'_0 \wedge dF_1 + dF_0 \wedge dF'_1 ).
$$

 Let $\eta \in \mathrm H^0(\p r, \Omega^1(d+2))$ represent
an element of $T_{\omega}\mathscr F$, that is, $d\omega
\wedge d\eta = 0$. We shall prove that $\eta$ belongs to the
image of $d \rho(F_0, F_1)$, \ie
$$ \eta = i_R ( dF'_0 \wedge dF_1 + dF_0 \wedge dF'_1 )
$$ for some $F'_0 \in \s {d_0}$ and $F'_1 \in \mathbf
S_{d_1}$.

Since $d \omega = dF_0 \wedge dF_1$, applying the division
Lemma \ref{L:2} to $d \eta$ it follows that there exist
homogeneous $1$-forms $\alpha$ and $\beta$ such that
$$ d\eta=\alpha \wedge dF_0 + \beta \wedge dF_1 .
$$ Notice that $d \eta$ is a 2-form with coefficients
homogeneous polynomials of degree $d = d_0 + d_1 - 2$. Hence
the coefficients of $\alpha$ (resp. $\beta$) are homogeneous
of degree $d_1 -1$ (resp.  $d_0 -1$). Applying exterior
derivative we find
$$ d\alpha \wedge dF_0 + d\beta \wedge dF_1 = 0 .
$$ Multiplying by $dF_1$ we get $d\alpha \wedge dF_0 \wedge
dF_1 = 0$.  From lemma \ref{L:2} applied to $d\alpha$ we
deduce
$$ d \alpha = \alpha' \wedge dF_0 + \alpha'' \wedge dF_1
$$ where $ \alpha'$ and $\alpha''$ are 1-forms with
coefficients homogeneous polynomials of respective degrees
$d_1 - 2 - (d_0 - 1) = d_1 - d_0 - 1$ and $d_1 - 2 - (d_1 -
1) = - 1$. Hence $\alpha'' = 0$. Similarly,
$$ d \beta = \beta' \wedge dF_0 + \beta'' \wedge dF_1
$$ where $\beta'$ and $\beta''$ are 1-forms with
coefficients homogeneous polynomials of respective degrees
$d_0 - 2 - (d_0 - 1) = - 1$ and $d_0 - 2 - (d_1 - 1) = d_0 -
d_1 - 1$. Hence $\beta' = 0$.

 Suppose that $d_0 = d_1$.  By the considerations above
regarding degrees, $\alpha' = \beta'' = 0$. Thus $\alpha$
and $\beta$ are closed $1$-forms. Therefore $\alpha = -
dF'_1$ and $\beta = dF'_0$ where $F'_i$ is some homogeneous
polynomial of degree $d_i$. It follows that $d\eta= dF'_0
\wedge dF_1 + dF_0 \wedge dF'_1 $ and since $i_R(d \eta) =
(d+1) \eta$ we obtain that $\eta$ is a scalar multiple of
$i_R ( dF'_0 \wedge dF_1 + dF_0 \wedge dF'_1 )$. Therefore
the Proposition is proved in the case $d_0 = d_1$.

 Now suppose $d_0 \ne d_1$, say $d_0 > d_1$. Then $d_1 - d_0
- 1 < 0$. Hence $d \alpha = 0$ and $d \beta = \beta'' \wedge
dF_1 $. Repeating the argument of the previous case we
obtain a sequence of $1$-forms $\beta_i$, $i \in \mathbb N$,
such that
$$ d \beta_i = \beta_{i+1} \wedge dF_1
$$ Comparing degrees it follows that, for $k \gg 0$,
$\beta_k=0$. Thus $ d \beta_{k-1}=0$ and there exists a
homogeneous polynomial $b_{k-1} $ such that $\beta_{k-1}=
db_{k-1}$. Then $d \beta_{k-2} =d b_{k-1} \wedge dF_1$ and
hence $\beta_{k-2} = b_{k-1} dF_1 + db_{k-2}$ for a suitable
homogeneous polynomial $b_{k-2}$. Then $d \beta_{k-3} =
\beta_{k-2} \wedge dF_1 = db_{k-2} \wedge dF_1$. Hence there
exists $b_{k-3}$ such that $\beta_{k-3} = b_{k-2} dF_1 +
db_{k-3}$. Iterating this, we conclude that $ \beta =
\beta_0 = b_1 dF_1 + db_0$ and therefore
$$ d\eta = dF'_1\wedge d F_0 + d F'_0 \wedge dF_1
$$ where $dF'_1 = \alpha$ and $dF'_0 = db_0 $, as wanted.
\end{proof}

\subsection{ Proof  of Theorem \ref{T:1}.}

\medskip

As a matter of fact we  prove  the
following slightly more precise statement.

\begin{thm}\label{T:reduzida}
If $r\ge 3$ then $\mathscr R(r, d_0, d_1)$ is an irreducible
component of $\mathscr F(r, d)$. Moreover, $\mathscr F(r,
d)$ is smooth and reduced at the points of $\mathcal U$.
\end{thm}
\begin{proof}
 Write as before $\rho: P \dashrightarrow\mathscr F$, where
$P = \ps{ \s {d_0} } \times \ps{ \s {d_1} }$,
$ \mathscr F = \mathscr F(r, d)$ and $\mathscr R = \mathscr
R(r, d_0, d_1)$ is the closure of the image of $\rho$. Put
$F=(F_0,F_1)\in P$. Proposition \ref{P:derivative} implies
that for $\omega = \rho(F)$, the derivative
$$ d \rho(F): T_FP \to T\mathscr F_{\omega}
$$ is surjective and also factors through $T_{\omega}
\mathscr R \subseteq T_\omega\mathscr F$. Then $T_\omega
\mathscr R = T_\omega\mathscr F$. It follows that $\mathscr
R$ is an irreducible component of $\mathscr F$ and $\mathscr
F$ is reduced at the generic point of $\mathscr R$.
\end{proof}

\section{Stability of quasi-homogeneous rational maps}\label{S:par}

In this section we exhibit some previously unknown
irreducible components $\mathscr R(r, d_0, \dots, d_q)$ of
$\mathscr F_q(r, d)$, generalizing the case $q=1$ of the
previous section.

A point of $\mathscr R(r, d_0, \dots, d_q)$ will be a
twisted $q$-form $\omega \in \mathrm H^0(\p r, \Omega^q(d+q+1))$
of type
\begin {equation} \label{omegaq}
\omega = i_R(dF_0 \wedge \dots \wedge dF_q) = \sum_{0 \le j
\le q} (-1)^j d_j F_j \ dF_0 \wedge \dots \wedge \widehat
{dF_j} \wedge \dots \wedge dF_q
\end {equation}
where $F_j \in \s {d_j}$ is a homogeneous polynomial
of degree $d_j$ in $r+1$ variables, and
\be\label{dqs}
d_0 +\dots + d_q = d + q + 1.
\ee
We call $\omega$ a {\it rational $q$-form} in $\p r$ of type
$(d_0, \dots, d_q)$.

More precisely, $\mathscr R(r, d_0, \dots, d_q)$ is defined
as the closure of the image of the rational map
\begin {equation}\label{mapq}
\rho: \ps{ \s {d_0} } \times \dots \times \ps{
\s {d_q} } \dashrightarrow \ps{ \mathrm H^0(\p r, \Omega^q(d+q+1)) }
\end {equation}
induced by the multilinear map
$$ \mu: \s {d_0} \times \dots \times \s {d_q}
\to \mathrm H^0(\p r, \Omega^q(d+q+1))
$$ such that $\mu(F_0, \dots, F_q) = i_R(dF_0 \wedge \dots
\wedge dF_q)$. The base locus of $\rho$ is described in
\rf{baselocusq} below.

As in the previous section, we define the open subset
\be\label{Uq} \mathcal U = \{ \omega \in \mathscr R(r, d_0,
\dots, d_q) \,|\,\mathrm{codim \ sing}(d\omega )\ge 3 \text{
and } \mathrm{codim \ sing}(\omega)\ge 2 \}. \ee

With notation as above, our main purpose in this section is to
prove the following Theorem \ref {T:reduzidaq},
which is a more precise version of  Theorem  \ref{T:2} of the Introduction.

\begin{thm} \label{T:reduzidaq} Suppose $r\ge 3$ and $1 \le q \le r-2$. Then
$\mathscr R(r, d_0, \dots, d_q)$ is an irreducible component
of $\mathscr F_q(r, d)$. Moreover, $\mathscr F_q(r, d)$ is
smooth and reduced at the points of $\mathcal U$.
\end{thm}

The strategy is the same as the one used to prove Theorem
\ref{T:reduzida}. Let us denote by $\mathscr F = \mathscr
F_q(r, d)$. The scheme $\mathscr F$ is defined by the
quadratic equations
\begin{equation} \label{schemeq}
i(v_J) \omega \wedge \omega = 0 \ \ \mathrm{and } \ \ i(v_J)
\omega \wedge d \omega = 0
\end{equation}
for all $J \subset \{0, \dots, r \}$ of cardinality $q-1$.

 The tangent space $T_\omega \mathscr F$ of $\mathscr F$ at
a point $\omega$ is represented by the forms $\omega' \in
\mathrm H^0(\p r, \Omega^q(d+q+1)) / (\omega) $ such that
$\omega_{\epsilon} = \omega + \epsilon \omega'$ satisfies
the conditions (\ref{schemeq}) modulo $\epsilon^2 $, that is
$$ i(v_J) \omega_{\epsilon} \wedge \omega_{\epsilon} = 0 \ \
\mathrm{and } \ \ i(v_J) \omega_{\epsilon} \wedge d
\omega_{\epsilon} = 0
$$ modulo $\epsilon^2 $, for all $J \subset \{0, \dots, r
\}$ of cardinality $q-1$. Expanding, one obtains
\begin{equation} \label{tangentq}
i(v_J) \omega' \wedge \omega + i(v_J) \omega \wedge \omega'
= 0 \ \ \mathrm{and } \ \ i(v_J) \omega' \wedge d\omega +
i(v_J) \omega \wedge d\omega' = 0 .
\end{equation}

In order to work out $\omega'$ from (\ref{tangentq}) we will
need a pair of technical results.

\subsection{Lemmata}

\medskip

The first technical Lemma is a generalization of Lemma
 \ref{L:2} that will be a central tool in the rest of this
 article.

\begin{lemma}l\label{L:saito+} Let $F_0 , \ldots, F_q$ be homogeneous {
polynomial} functions on $\mathbb C^{r+1}$ and let $\Theta
\in \Omega^{q+1}(\mathbb C^{r+1})$ be the $(q+1)$-form given
by
\[
\Theta = dF_0 \wedge \ldots \wedge dF_q \, .
\]
Suppose that $\mathrm{codim } \, \mathrm{ sing} (\Theta) \ge
3$. If $\eta \in \Omega^{q+1}(\mathbb C^{r+1})$ is such that
$ \eta \wedge dF_i \wedge dF_j =0$ for every $0\le i < j \le
q$ then there exist holomorphic $1$-forms $\alpha_0, \ldots,
\alpha_q \in \Omega^1(\mathbb C^{r+1})$ such that
$$ \eta = \sum_{i=0}^q \alpha_i \wedge dF_0 \wedge \ldots
\widehat{dF_i} \ldots \wedge dF_q .
$$
\end{lemma}
\begin{proof}
For the second item let $\mathcal U$ be an open covering of
$\mathbb C^{r+1} \setminus \mathrm{sing}(\Theta)$. Since
$\mathrm{codim } \, \mathrm{ sing} (\Theta) \ge 3$ we can
assume that over each open set $U \in \mathcal U$ our set of
functions is part of a coordinate system on $U$. It is then
clear that
\[
\eta _{|U} = \sum \alpha_{i,U} \wedge dF_0 \wedge \ldots
\widehat{dF_i} \ldots \wedge dF_q
\]
for suitable $1$-forms $\alpha_{0,U},\ldots, \alpha_{q,U}
\in \Omega^1 (U)$.

A simple computation shows that over $U\cap V$
\[
(\alpha_{i,U} - \alpha_{i,V}) \wedge \Theta =0 \, .
\]
It follows from Saito's Lemma \cite{Saito} that there exists
a unique $(q+1) \times (q+1)$ matrix $A_{U\cap V}$ with
entries in $\mathcal O(U\cap V)$ such that
\[
\left[\begin{array}{c} \alpha_{0,U} - \alpha_{0,V} \\ \vdots
    \\ \alpha_{q,U} - \alpha_{q,V}
\end{array}
\right] = A_{U\cap V} \cdot \left[\begin{array}{c} dF_0 \\
\vdots \\ dF_q \\\end{array} \right]
\]
Of course the collection of matrices $A_{U\cap V}$ with
$(U,V)$ ranging in $\mathcal U^2$ defines an element of
$\mathrm H^1( \mathbb C^{r+1} \setminus
\mathrm{sing}(\Theta) , \mathbb M \otimes \mathcal O) \cong
\mathrm H^1( \mathbb C^{r+1} \setminus \mathrm{sing}(\Theta)
, \mathcal O) \otimes \mathbb M$, with $\mathbb M$ being the
vector space of $(q+1)\times (q+1)$ matrices.

The hypothesis $\mathrm{codim} \ \mathrm{sing}(\Theta) \ge
3$ implies that this cohomology group is trivial, see for
instance \cite[pg. 133]{Grauert}. Therefore we may write
$A_{U\cap V} = A_U - A_V$ where $A_U,A_V$ are matrices of
holomorphic functions in $U$ resp. $V$. We can thus set
\[
\left[\begin{array}{c} \alpha_{0} \\ \vdots \\ \alpha_{q} \\
\end{array}
\right] = \left[ \begin{array}{c} \alpha_{0,U} \\ \vdots \\
\alpha_{q,U} \\\end{array} \right] - A_U \cdot
\left[\begin{array}{c} dF_0 \\ \vdots \\ dF_q \\\end{array}
\right] = \left[ \begin{array}{c} \alpha_{0,V} \\ \vdots \\
\alpha_{q,V} \\\end{array} \right] - A_V \cdot
\left[\begin{array}{c} dF_0 \\ \vdots \\ dF_q \\\end{array}
\right]
\]
as the sought global $1$-forms at least over $\mathbb
C^{r+1} \setminus \mathrm{sing}(\Theta)$. To conclude one
has just to invoke Hartog's extension Theorem to ensure that
these $1$-forms extend to $\mathbb C^{r+1}$.
\end{proof}

By expanding in its homogeneous components both sides of the
equality
$$
\eta = \sum_{i=0}^q \alpha_i \wedge dF_0 \wedge \ldots
\widehat{dF_i} \ldots \wedge dF_q .
$$
it can be easily seen that if $\eta$ is a homogeneous
polynomial $q$-form then the $1$-forms $\alpha_0, \ldots,
\alpha_q$ can be assumed homogeneous polynomial $1$-forms.

\medskip

The second technical Lemma in this subsection replaces the
iteration argument in the proof of Theorem \ref{T:reduzida}

\begin{lemma} \label{lemmainductivo} For $j= 0, \dots, q$ let $F_j \in
\s {d_j}$ be a homogeneous polynomial of degree
$d_j$.  Suppose $\omega = i_R(dF_0 \wedge \dots \wedge
dF_q)$ satisfies $\mathrm{codim \ sing \, }(d\omega ) \ge
3$. Then, for $\alpha \in \mathrm H^0(\p r, \Omega^1(e))$
the following conditions are equivalent:
\begin{enumerate}[\rm(a)]
\item
$d\alpha = \sum_{0 \le k \le q} A_k \wedge dF_k$ for some
$A_k \in \mathrm H^0(\p r, \Omega^1(e - d_k))$.
\item
$\alpha = dG + \sum_{0 \le k \le q} H_k \ dF_k$ for some $G
\in \s e$ and $H_k \in \s {e - d_k}$.
\end{enumerate}
\end{lemma}

\begin{proof}
It is clear that (b) implies (a). Let us prove the converse,
by induction on $e \in \mathbb N$. If (a) holds, applying
exterior derivative we get
$$ 0 = d^2 \alpha = \sum_{0 \le k \le q} d A_k \wedge dF_k
\implies dA_k \wedge dF_0 \wedge \cdots \wedge d F_q =0 .
$$ By the hypothesis on the $F_j$ and Lemma \ref{L:2},
$$ d A_k = \sum_{0 \le h \le q} A_{kh} \wedge dF_h
$$ for some $A_{kh} \in \mathrm H^0(\p r, \Omega^1(e - d_k -
d_h))$.  Since $e - d_k < e$, the inductive hypothesis
applies to $A_k$ and yields
$$ A_k = dG_k + \sum_{0 \le h \le q} H_{kh} \ dF_h
$$ for some $G_k \in \s {e - d_k}$ and $H_k \in
\s {e - d_k - d_h}$.  Replacing in (a) we find
$$ d\alpha = \sum_{k} dG_k \wedge dF_k + \sum_{h, k} H_{kh}
\ dF_h \wedge dF_k.
$$ Since $i_R \alpha = 0$, we have \,$e\cdot\alpha = i_R d
\alpha$.  Applying $i_R$ we obtain, after a little
calculation
$$ e \cdot \alpha = dG + \sum_{0 \le k \le q} H_k \ dF_k
$$ where
$$ G = - \sum_{k} d_k F_k G_k, \ \ \ H_k = (d_k + e) G_k +
\sum_{h} d_h F_h (H_{kh} - H_{hk})
$$ as claimed.
\end{proof}

\subsection{Surjectivity of the derivative and  proof of Theorem  \ref{T:2}}

\medskip
Now we are ready to complete the proof of Theorem \ref {T:reduzidaq}
and hence of  Theorem \ref{T:2} of the Introduction.
The proof  follows from  Proposition \ref{P:derivativeq}
below combined with the same argument used in the proof of Theorem \ref{T:reduzida}.

\begin{prop} \label{P:derivativeq}
Suppose $r\ge 3$ and  $1 \le q < r-1$.
If $\u{F}=(F_0, \dots, F_q) \in \prod_i \ps{ \s {d_i} }$ is
such that $\rho(\u{F}) = \omega \in \mathcal U$ then the derivative
$$
d \rho (\u{F}) : T_{\u{F}}(\ps{ \mathbf S_{d_0} }\times \cdots \times \ps{ \s {d_q} }) \to
T_{\omega} \mathscr F
$$
is surjective.
\end{prop}
\begin{proof}
At a point $\u{F}=(F_0, \dots, F_q)$ belonging to the domain
of $\rho$ the derivative
\begin {equation} \label{derivativeq}
d \rho (\u{F}): \s {d_0} / (F_0) \times \dots \times
\s {d_q} / (F_q) \rightarrow T_\omega \mathscr F
\end {equation}
is calculated by multilinearity as
$$ d \rho(\u{F}) (F'_0, \dots, F'_q) = \sum_{0 \le j \le q}
i_R (dF_0 \wedge \dots \wedge dF'_j \wedge \dots \wedge
dF_q).
$$

Let $\omega = \rho(\u{F}) \in \mathcal U$ and $\omega' \in
T_\omega \mathscr F$. From \rf{tangentq} we have
$$ i(v_J) \omega' \wedge d\omega = - i(v_J) \omega \wedge
d\omega' .
$$ Since $d \omega$ is a constant multiple of $dF_0 \wedge
\dots \wedge dF_q$ (see Lemma \ref{L:1} ), by exterior
multiplication with $dF_j$ we obtain
$$ dF_j \wedge i(v_J) \omega \wedge d\omega' = 0
$$ for all $j, J$.

Let $Y_j, (0 \le j \le q)$, be rational vector fields such
that $dF_i(Y_j) = \delta_{ij}$. For $J = \{0, \dots, q\}
\setminus \{i, j\}$ we have $ i(v_J) \omega = \lambda (F_i
dF_j - F_j dF_i)$. Then, \\\centerline{ $0 = dF_j
\wedge i(v_J) \omega \wedge d\omega' = \lambda dF_j \wedge
F_j dF_i \wedge d\omega'$,} \\ which implies that
$$ dF_i \wedge dF_j \wedge d\omega' = 0
$$ for all $0 \le i, j \le q$.

Lemma \ref{L:saito+} implies that
\begin {equation} \label{domegaprime}
d\omega' = \sum_{0 \le j \le q} \alpha_j \wedge dF_0
\wedge\dots \wedge \widehat {dF_j} \wedge \dots \wedge dF_q
\end {equation}
for some $\alpha_j \in \mathrm H^0(\p r,
\Omega^1(d_j))$. Applying exterior derivative we find
$$ 0 = d^2\omega' = \sum_{0 \le j \le q} d\alpha_j \wedge
dF_0 \wedge \dots \wedge \widehat {dF_j} \wedge \dots
\wedge{}dF_q.
$$ Taking wedge product with $dF_j$ we get
$$ d\alpha_j \wedge (dF_0 \wedge \dots \wedge dF_q) = 0
$$ for all $j$. Therefore, thanks to Lemma \ref{L:2},
$$ d\alpha_j = \sum_{0 \le k \le q} A_{jk} \wedge dF_k
$$ for suitable $A_{jk} \in \mathrm H^0(\p r, \Omega^1(d_j -
d_k))$.  Lemma \ref{lemmainductivo} implies that
$$ \alpha_j = dG_j + \sum_{0 \le k \le q} H_{jk} \ dF_k
$$ for some $G_j \in \s {d_j}$ and $H_{jk} \in \s {d_j - d_k}$
(we use the convention $\s e = 0$ for $e < 0$).
Replacing in \rf{domegaprime} above we have
\begin {equation} \label{domegaprime2}
d\omega' = \sum_{0 \le j \le q} dG_j \wedge dF_0 \wedge
\dots \wedge \widehat {dF_j} \wedge \dots \wedge dF_q +
c \ dF_0 \wedge \dots \wedge dF_q
\end {equation}
for some $c \in \mathbb C$. Since $i_R \omega' = 0$, Lemma
\ref{L:1} yields $(\sum_{i} d_i) \ \omega' = i_R d
\omega'$. Applying $i_R$ to \rf{domegaprime2} and taking
\rf{derivativeq} into account, we obtain
$$ \omega' = d \rho(\u{F}) (F'_0, \dots, F'_q)
$$ where $F'_j = \frac{ (-1)^j}{(\sum_{i}
d_i)\vphantom{I^I}} \ G_j$.  Therefore $d \rho (\u{F})$ is
surjective, as claimed. \end{proof}

\section{Geometry of the parametrization}
In this section we  analyze the parametrization
$$ \rho: \ps{ \s {d_0} } \times \dots \times \ps{
\s {d_q} } \dashrightarrow \mathscr R_q(r, \bar d)
\subset \ps{ \mathrm H^0(\p r, \Omega^q(d+q+1)) } \, ,
$$ where $\s {d_i}= \mathrm H^0 (\mathbb P^r,
\mathcal O_{\mathbb P^r} (d_i))$, $d = \sum d_i$ and $\bar d
= (d_0 , \ldots, d_q)$.

\subsection{Base locus}

\medskip

Let us start by describing the base locus $\mathbf B (\rho)$
of $\rho$.

If $i_R(dF_0 \wedge \dots \wedge dF_q) = 0$,
applying exterior differentiation and Lemma\,\ref{L:1} we
obtain that $dF_0 \wedge \dots \wedge dF_q = 0$. This means
that the Jacobian matrix of $F_0, \dots, F_q$ has rank $<
q+1$ everywhere, that is, the derivative of the map
$$
\u F: \mathbb C ^{r+1} \to \mathbb C ^{q+1}
$$ defined by $\u F(x) = (F_0(x), \dots, F_q(x))$ has rank $<
q+1$ at every $x \in \mathbb C ^{r+1}$. This is equivalent
to the fact that $F$ is not dominant, that is, $f(F_0,
\dots, F_q) = 0$ for some non-zero polynomial $f \in \mathbb
C[y_0, \dots, y_q]$ (\ie the $F_j$ are algebraically
dependent). We thus obtain
\begin{equation}\label{baselocusq}
\mathbf B (\rho) =
\{(F_0, \dots, F_q)
\in \prod_i \ps{ \s {d_i} }\,|\, \u F: \mathbb C ^{r+1} \to
\mathbb C ^{q+1} \ \mathrm{is \ not \ dominant} \}.
\end{equation}

For $q=1$ the set theoretical description of $\rho$ is
rather simple:
\begin{equation} \label{baselocus}
\mathbf B (\rho) = \{(F_0, F_1) \in \ps{ \mathbf{S}_{d_0} }
\times \ps{ \s {d_1} } \,|\, F_0^{d_1} = F_1^{d_0} \}
\, .
\end{equation}

For general $q$ we have a stratification
$$ \mathbf B (\rho)_1 \subset \mathbf B (\rho)_2 \subset
\dots \subset \mathbf B (\rho)_q =\mathbf B (\rho)
$$ where $\mathbf B (\rho)_k = \{(F_0, \dots, F_q) \,|\,
\mathrm{\dim{} image}(F) \le k \}$. The first stratum
$\mathbf B (\rho)_1$ is set-theoretically equal to
\[
 \{ (F_0, \ldots, F_q) \in\prod_i  \ps{ \s {d_i} } \, |
\, F_0^{\hat d_0} = \ldots = F_q ^{\hat d_q} \}
\]
where $\hat d_j = \prod_{i\ne j} d_i$. For $k>1$ the same
set theoretical description is considerably more complex and
we will carry it out only in very particular cases in \S
\ref{S:section Degree}.

Beware that the scheme structure of $\mathbf B(\rho)$ is
often non-reduced, see \S \ref{SS:23}.

At any rate, we register the following easy consequence of
Lemma\,\ref{L:1}.
\begin{prop}\label{P:rorotil}
Let
$$\ba{rcc}
\wt\rho:\prod_i \ps{\s{d_i}} &\dasharrow&
\ps{\s{d-1}\otimes\wed{q+1}\s1\us}\\
(F_0,\dots,F_q)&\mapsto&dF_0\wedge\dots\wedge F_q.
\ea$$
Then the base loci of \wt\rho\ and $\rho$ are one and the
same as schemes.
\end{prop}

\begin{proof}
Let $V\subset\s{e}\otimes\wed q\s1\us$ be the subspace of
closed $q$--forms with coefficients of degree $e$. Put
$W=i_R(V)\subset\s{e+1}\otimes\wed{q-1}\s1\us$. Then
$i_R:V\ra W$ is a linear isomorphism. We still denote by
$i_R:\ps V\ra \ps W$ the projectivization. Since the image
of \wt\rho\ lies in \ps V and
$\rho=i_R\circ\wt\rho$, the assertion follows.
\end{proof}

\subsection{Weighted homogeneous polynomials}\label{S:wp}

\medskip

  Fix $\bar d = (d_0, \dots, d_q) \in \mathbb N^{q+1}$ and
$e \in \mathbb N$. A polynomial $f$ in $\mathbb C[{y_0,
\dots, y_q}]$ is said to be weighted homogeneous of type
$\bar d$ and degree $e$ if
$$f(\lambda^{d_0} y_0, \dots, \lambda^{d_q} y_q) =
\lambda^{e} f(y_0, \dots, y_q)
$$ for any $\lambda \in \mathbb C$. Equivalently, $f$ is a
linear combination of monomials
$$\prod_{0 \le j \le q} y_j^{\alpha_j} \mathrm{\ \ such \
that \ \ } \bar d\cdot\aa:=\sum_{0 \le{} j\le q} d_j
\alpha_j = e.
$$
This is tantamount to declaring each variable $y_i$ to be of
degree $d_i$.

We denote by $$ \s {q,\bar d, e} $$ the $\mathbb
C$-vector space of all such polynomials and write its
dimension as $N(q, \bar d, e)$. Notice that $N(q,\bar
d,e)=\dim\s {q,\bar d,e}$ can be expressed by the
Hilbert series
$$ H(t)=\sum_e N(q,\bar d,e)t^e=\frac
1{\vphantom{I^{I^I}}\prod_{i=1}^{q} (1-t^{d_i})}.
$$

Throughout we will assume that the vector of natural numbers
$\bar d \in \mathbb N^{q+1}$ is non-decreasingly ordered,
\ie $d_0\leq d_{1}\leq \cdots\leq d_q$.

Define $\bar e = \bar e(\bar d) = (e_1 , \ldots , e_k)$ such
that $e_i < e_{i+1}$ and $\cup_{0\le i\le q} \{d_i\} =
\cup_{1\le i\le k} \{e_i\}$. If $n_i$ stands for the number
of times the natural number $e_i$ appears in $\bar d$ then the
pair $(\bar e, \bar n )$, where $\bar n = (n_1,\ldots
,n_k)$, determines $\bar d$.

Set $ q_j = -1 + \sum_{1\le i \le j} n_i$, and for $l = 1,  \dots, k$
\[
 \quad \bar d_l =
  (\underbrace{e_1, \ldots, e_1}_{ n_1 \text{ times}},
  \underbrace{e_2, \ldots, e_2}_{ n_2 \text{ times}},
  \ldots, \underbrace{e_l, \ldots, e_l}_{ n_l \text{ times}}).
\]
Clearly, for each $f \in \s {q,\bar d, e_j}$, no
variable $y_i$ with weight $d_i>e_j$ occurs in $f$; thus
$$ \displaystyle{\s {q,\bar d, e_j} \cong \mathbf
S_{q_j , \bar d_j , e_j }}.
$$

Denote by $\mathbb E^{q+1}= \mathrm{End}(\mathbb C^{q+1})$
the set of all polynomial maps $f: \C^{q+1} \to
\C^{q+1}$. It is a ring under
sum and composition of maps. If $f = (f_0, \dots, f_q) \in
\mathbb E^{q+1}$, we say that $f$ is of type $\bar d$ if
$f_i$ is weighted homogeneous of type $\bar d$ and degree
$d_i$, for all $i=0, \dots, q$.

\begin{lemma}\label{F:1} Maps of type $\bar d$  form a
subring of $\mathbb E^{q+1}$. More precisely, if $f, g \in
\mathbb E^{q+1}$
are of type $\bar d$ then $f \circ g$ is of type $\bar d$.
Moreover, the set
$$\mathrm{GL}(q,\bar d) = \{f \in \mathbb E^{q+1} | f
\text{ is of type } \bar d \text{ and } df(0) \text{ is
invertible} \}$$ is a group.
\end{lemma}

\begin{proof}
$(f_i \circ g) (t^{d_0} y_0, \dots, t^{d_q} y_q) = f_i (g_0
(t^{d_0} y_0, \dots, t^{d_q} y_q)), \dots, g_q (t^{d_0} y_0,
\dots, t^{d_q} y_q)) = f_i (t^{d_0} g_0(y_0, \dots, y_q),
\dots, t^{d_q} g_q(y_0, \dots, y_q)) = t^{d_i} f_i (g_0(y_0,
\dots, y_q), \dots, g_q(y_0, \dots, y_q)) = t^{d_i} (f_i
\circ g) (y_0, \dots, y_q),$

We have  $\mathrm G=\mathrm{GL}(q,\bar
d)$ is closed under compositions. It remains to show that
every element is invertible in $\mathrm G$.  Let us denote
the block of variables of weight $e_i$ by
$$
\u y_{1}=\underbrace{
y_{0},\dots,y_{q_1}}_{\text{ (weight }e_1)},\quad
\u y_{2}=\underbrace{
y_{q_1+1},\dots,y_{q_2}}_{\text{ (weight }e_2)},
\quad\dots,\quad
\u y_{k}=\underbrace{
y_{q_{k-1}},\dots,y_{q_k}}_{\text{ (weight }e_k)}.
$$
The main point is that each $f \in \mathrm G$ has the
following triangular shape,
$$
(\u f_1(\u y_1),\u f_2(\u y_1,\u y_2),\dots,
\u f_{k}(\u y_1,\dots,\u y _k)).
$$
Here
$$
\u f_i(\u y_1,\dots,\u y_i)=(f_{i1}(\u y_1,\dots,\u y_2),
\dots,f_{2n_i}(\u y_1,\dots,\u y_i)),
$$
with
$$
f_{ij}(\u y_1,\dots,\u y_i)=g_{ij}(\u y_1,\dots,\u y_{i-1})
+h_{ij}(\u y_i)\in\s {q_i,\bar d_i,e_i}
$$
where
$h_{ij}(\u y_i)$ is in fact linear in the block of
variables $\u y_i$ of weight $e_i$.
Indeed, since $e_{i+1}>e_i$, no $\u y_{i+1}$ occurs in
$\u f_{i}$.
Thus $f$ can be written as
$$
(\u h_1(\u y_1),\u h_2(\u y_2)+\u g_2(\u y_1),\dots,
\u h_k(\u y_k)+\u g_{k}(\u y_1,\dots,\u y _{k-1})).
$$
Now we see that $df(0)$ is made up of blocks of the linear
maps $\u h_i=d\u h_{i}:\C^{n_i}\ra\C^{n_i}$. Hence
invertibility of the former is equivalent to $d\u
h_{i}\in\mathrm{GL}_{n_i}\,\forall i$. Thus, given
$(z_1,\dots,z_q)=(f(y))$, one can solve successively
$$
\left\{\ba l\u y _1= \u h_{1}\inv(\u z_1), \text{ then }\\
\u y_2=\u h_2\inv(\u z_2-\u g_2(\u y_1)),\\\vdots\\
\u y_k=\u h_k\inv(\u z_k-\u g_k(\u y_1,\dots,\u y _{k-1})).
\ea\right.
$$
\end{proof}

The group $\mathrm{GL(q,\bar d)}$ naturally acts on
the domain of $\mu$ (cf.\,\ref{mapq}):
\begin{eqnarray*}
\mathrm{GL(q,\bar d)} \times \prod_{0 \le j \le q} { \mathbf
S_{ d_j}} &\longrightarrow & \prod_{0 \le j \le q} { \mathbf
S_{d_j}} \\ \left( f , (F_0, \ldots, F_q) \right) &\mapsto&
(f_0 (\u F) , \ldots, f_q(\u F)) \, .
\end{eqnarray*}
In other words, considering $\u F$ as a polynomial map
$\u F: \mathbb C^{r+1} \to \mathbb C^{q+1}$, the action
is just composition with a polynomial map
$f: \mathbb C^{q+1} \to \mathbb C^{q+1}$ which belongs
to $\mathrm{GL(q,\bar d)}$.

\subsection{The fibers of $\rho$}

\medskip

The key tool for the description of the fiber of $\rho$ and
the proof of Theorem \ref{T:rational} is the following
Proposition.

\begin{prop} \label{P:fiberq} Let $\u F = (F_0, \dots, F_q), \u G= (G_0,
\dots, G_q) \in \s {d_0} \times \dots \times \mathbf
S_{d_q} $ Suppose that both $dF_0 \wedge \cdots \wedge dF_q$
and $dG_0 \wedge \cdots \wedge dG_q$ are non-zero
$(q+1)$-forms. If $\cod \mathrm{sing}(dF_0 \wedge \cdots
\wedge dF_q) \ge 2$ then the following conditions are
equivalent:
\begin{enumerate}[{\rm(a)}]
\item
{$i_R(dF_0 \wedge \dots \wedge dF_q) = i_R(dG_0 \wedge \dots
\wedge dG_q)$ up to a constant multiple.}

\item
{$dF_0 \wedge \dots \wedge dF_q = dG_0 \wedge \dots \wedge
dG_q$ up to a constant multiple.}

\item
$dG_j = \sum_{0 \le k \le q} A_{jk} \ dF_k$ for some $A_{jk}
\in \s {d_j - d_k}$, for all $j$.

\item
$G_j = f_j(F_0, \dots, F_q)$ for some $f_j \in \mathbb
C[y_0, \dots, y_q]$, for all $j$.

\item
$G_j = f_j(F_0, \dots, F_q)$, for all $j$ for a unique $f_j
  \in \s {q,\bar d,d_j}$.  Moreover,
  $(f_0,\dots,f_q)$ belongs to $\mathrm{GL}(q,\bar d)$.
\end{enumerate}
\end{prop}
\begin{proof}
$\rm(a) \Leftrightarrow (b)$: Use the identity $d(i_R(dF_0
\wedge \dots \wedge dF_q))= (q+d)(dF_0 \wedge \dots \wedge
dF_q) $ from Lemma\,\ref{L:1}.

$\rm(b) \Rightarrow (c)$: Multiplying by $dG_j$ we obtain
$dG_j \wedge dF_0 \wedge \dots \wedge dF_q = 0$. Since $\u
F$ is generic, it follows by the division lemma that the
$dG_j$ are linear combinations of the $dF_k$. The
coefficients may be chosen as homogeneous polynomials,
necessarily of the stated degree.

$\rm(c) \Rightarrow (b)$: Using the hypothesis and
calculating wedges we have
$$ dG_0 \wedge \dots \wedge dG_q = \mathrm{det}(A) \ dF_0
\wedge \dots \wedge dF_q.
$$
Now $\mathrm{det}(A)$ is a non-zero homogeneous polynomial,
and its degree is zero, so it is a constant, thereby proving
the claim.

$\rm(d) \Rightarrow (e)$: Let $f_j = \sum_{\alpha}
c_{\alpha} y^{\alpha}$, where $\alpha \in \mathbb N^{q+1}$
and $c_{\alpha} \in \mathbb C$, so that $G_j = \sum_{\alpha}
c_{\alpha} F^{\alpha}$.  Write $f_j=g_j+h_j$ where $g_j$ is
the sum over the exponents \aa\ such that $\bar
d\cdot\aa=d_j$. We have $h_j(\u F)=0$ by the homogeneity of
$G_j$ and of the $F_k$. Therefore we may take $f_j=g_j$, the
weighted homogeneous polynomial that we needed.  Uniqueness
is clear since the $F_k$ are algebraically
independent. Finally, setting $f=(f_0,\dots,f_q)$, since
$$dG_0\wedge\dots\wedge{}dG_q=\det(df)dF_0\wedge\dots
\wedge{}dF_q
$$ it follows that $\det(df)=\det(df(0))$ is a nonzero
constant.

$\rm(e) \Rightarrow (d)$: obvious.

$\rm(d) \Rightarrow (c)$: If $G_j = \sum_{\alpha} c_{\alpha}
F^{\alpha}$, taking exterior derivative we immediately get
$dG_j$ as a linear combination of the $dF_k$.

$\rm(c) \Rightarrow (d)$: {It suffices to use Lemma
\ref{again} below. }
\end{proof}

\begin{lemma} \label{again} Let $\u F = (F_0, \dots, F_q) \in \mathbf
S_{d_0} \times \dots \times \s {d_q} $ be
generic. Let G be a homogeneous polynomial of degree e such
that $dG = \sum_{0 \le k \le q} A_{k} \, dF_k$ for some
$A_{k} \in \s {e - d_k}$. Then $G = f(F_0, \dots,
F_q)$ for a unique polynomial $f \in \s {q,\bar
d,e}$.
\end{lemma}
\begin{proof}
We proceed by induction on $e$. The assertion is clear for
$e=0$.  Taking exterior derivative we have $d^2G=\sum_k
dA_{k} \wedge dF_k = 0$. Thus $dA_{k} \wedge dF_0 \wedge
\dots \wedge dF_q=0$ for all $k$. Since $\u F$ is generic,
we get $dA_{k} = \sum_h B_{kh}\, dF_h$ for some $B_{kh} \in
\s {e - d_k - d_h}$. By the inductive hypothesis,
$A_{k} = f_{k}(F_0, \dots, F_q)$ for some polynomial
$f_{k}$. On the other hand, applying $i_R$ to $dG = \sum_k
A_{k} \, dF_k$ we obtain $e G = \sum_k A_{k}\, d_k
F_k$. Replacing here $A_{k}$ by $f_{k}(F_0, \dots, F_q)$ we
obtain the claim. Uniqueness and weighted homogeneity were
argued before.
\end{proof}

\begin{prop}\label{P:fibra}  For general $\u F = (F_0, \dots, F_q) \in
\prod_{0 \le j \le q} \s {d_j}$ we have a bijective
map
$$ \ba {ccc} \mathrm{GL}(q,\bar d)&\xymatrix{\ar@{->}[r]&} &
\mu^{-1} \mu (\u F) \\ (f_0,\dots,
f_q)&\xymatrix{\ar@{|->}[r]&}&(f_0(\u F), \dots,f_q(\u
F))\ea
$$ with $\mu$ the multilinear map inducing $\rho$ as in
\rf{mapq}.
\end{prop}
\begin{proof}
The assertion follows from the equivalence (a)$\iff$(e) in
\ref{P:fiberq}.
\end{proof}

\begin{cor}\label{fibdim} We have the formula for the fiber dimension,
$$ \dim\rho^{-1} \rho (\u F) = \sum_{0 \le j \leq q} (N(q,
\bar d, d_j) - 1).
$$
\end{cor}

\subsection{A natural factorization and proof of Theorem
\ref{T:rational}}\label{S:X}

\medskip

We will now proceed to describe a tower of open subsets of
Grassmann bundles birational to $\mathscr R(r, \bar d)$. We
preserve the notation of Subsection\,\ref{S:wp}.

Start with $Y_0 = G(n_1, \s {e_1})$, the grassmannian
of $n_1$-planes in $\s {e_1}$.  Let $X_1 \subset Y_1$
be the open subset defined as
\[
X_1 = \{ F_1 \wedge \cdots \wedge F_{n_1} \in G(n_1,
S_{e_1}) \, | \cod \mathrm{sing}(dF_0 \wedge \cdots
\wedge dF_{n_1}) \ge 2
\}.
\]
Now let $\cl A_2\ra X_1$ be the vector subbundle of the
trivial bundle $\s {e_2}\vez X_1$ with fiber over $\u
F_1=F_1 \wedge \cdots \wedge F_{n_1} \in X_1$ given by
$$
\cl A_2(\u F_1)=\{G\in\s {e_2}\st
dF_1 \wedge \cdots \wedge dF_{n_1} \wedge dG=0\}.
$$
Recalling Lemma\,\ref{L:2}(a), and the above considerations
on weighted homogeneity, we have in fact
$$
\cl A_2(\u F_1)=\{G\in\s {e_2}\st
G=f(\u F_1),f\in\s {q_1,\bar d_1, e_2}
\}\cong\s {q_1,\bar d_1, e_2}.
$$

Let $Y_2=G(n_2,\s {e_2}/\cl A_2)$
be the Grassmann bundle over $X_1$.
Notice that, for an element $\u
G_2= [G_1]\wedge \cdots \wedge [G_{n_2}] \in G(n_2,\mathbf
S_{e_2}/ \s {q,\bar d, e_2}(p))$ over a point $\u
F_1=F_1 \wedge \cdots \wedge F_{n_1} \in X_1$, the
$(n_1+n_2)$-form
\[
\eta(\u G_2) = dF_1 \wedge \cdots dF_{n_1} \wedge dG_1
\wedge \cdots \wedge dG_{n_2}
\]
is well-defined up to a non zero multiplicative
constant. Therefore we can set $X_2 \subset Y_2$ as the open
subset defined by
\[
X_2 = \{ \u G_2 \in Y_2 \, | \, \mathrm{codim \ sing \, } \eta(\u G_2)
\geq 2 \}
\]
Continuing, we have a vector subbundle $\cl A_3$ of $\mathbf
 S_{e_3}\times X_2$ with fiber
$$
\cl A_3(\u F_1,\u G_2)=\{H\in\s {e_3}\st
dF_1 \wedge \cdots \wedge dF_{n_1} \wedge
dG_1 \wedge \cdots \wedge dG_{n_2} \wedge dH=0\}.
$$
As before, this is isomorphic to $\s {q_2,\bar d_2,
  e_3}$ . Proceeding this way, we arrive at an open subset
$X=X_k\subset Y_k$ where $Y_k\ra X_{k-1}$ is the Grassmann
bundle  $G(n_k,\s {e_k}/\cl A_{k-1})$.  Clearly
$X$ is a rational variety just like all
Grassmann bundles over rational varieties. Using Proposition
\,\ref{P:fiberq}, we arrive at a birrational map
from  $X$ to $\mathscr R(r, \bar d)$.
It follows that $\mathscr R(r, \bar d)$ is rational and this
concludes the proof of Theorem \ref{T:rational} \qed

\section{Degree calculations} \label{S:section Degree}
Let $\bar d = (d_0,\ldots, d_q)$, $\bar e$, $\bar n$, \ldots
be as in the previous section. Here we proceed to find the
degree of the projective variety
$$
\mathscr R(r, \bar d) \subset \ps{ \mathrm H^0(\p r,
\Omega^q(d+q+1)) }
$$ in some cases.
We shall time and again profit from the following
consequence of Proposition\,\ref{P:rorotil}.
We consider
$$
\wt\rho: \prod_i  \ps{\s{d_i}}\dasharrow
\wt{\mathscr R}(r, \bar d)=(i_R)\inv\mathscr R(r, \bar d)
\subset
\ps{\s{d-1}\otimes\wed{q+1}\s1\us}.
$$
Thus we see that all degree calculations can be lifted from
$\ps{W} \subset \ps{\mathrm H^0(\p r,\Omega^q(d+q+1))
}$ to \ps V.

\subsection{Linear projections of grassmannians}

\medskip

When $q_1 = q$, i.e.  all the degrees $d_i$ are equal to
$e_1$, the variety $X$ constructed in \S \ref{S:X} is an
open subset of the grassmannian $G(q,\s {e_1})$. It
follows that the morphism $\bar \rho: X \to \mathscr R(r,
\bar d)$ gives rise to a rational map
\[
\wt \rho: G(q+1,\s {e_1}) \dashrightarrow \wt{\mathscr
R}(r, \bar d) \subset \ps{\s{d-1}\otimes\wed{q+1}\s1\us}.
\]
Notice that $\bar \rho$ is the
composition of Pl\"{u}cker's embedding with a central
projection
$$\ba{ccc}\ps{\bigwedge^{q+1} S_{e_1}} &\dashrightarrow&
\ps{\s{d-1}\otimes\wed{q+1}\s1\us}\\\na2
F_0\wedge\cdots\wedge F_q & \mapsto &
dF_0\wedge\cdots\wedge dF_q.
\ea$$

It is a simple exercise to show that $G(q+1,\mathbf
S_{e_1})$ is disjoint from the center of this projection if,
and only if, $q=1$ or $d_0=\cdots=d_q=1$. In both cases the
degree of these components is equal to the degree of the
corresponding grassmannians under Pl\"{u}cker's
embedding (see e. g. \cite{Kleiman}). More precisely, setting $N=(q+1)(r-q)=\dim
G(q+1,r+1)$, we have
\be\label{q1..1}
\boxed{
\begin{array}{lclcl}
     \deg( \mathscr R(q,1, \ldots,1) ) &=&
\deg G(q+1,\s {1})&=&
\frac{1!2!\cdots{}q!N!}{(r-q)!(r-q+1)!\dots r!}
\\\na5 \deg( \mathscr
     R(1,d,d) ) &=&\deg G(2,\s {d_1})
 &=& \frac1{N_d-1}\binom{2N_d-2}{N_d},\\\na3
\text{where}&&N_d=\binom{r+d}r-1.&&
\end{array}}
\ee

\medskip

\begin{remark}\label{blocus}\em
The scheme-theoretic structure of the base
locus of a rational map $\phi:Y\rar\ps{\C^N}$
is defined as follows (cf.\,\cite[7.17.3,\,p.\,168]{Hartshorne}).  We
are given a line bundle (=invertible sheaf) \cl L
over $Y$ together with a homomorphism
$\O_Y^{N}\ra \cl L$, surjective over the open
dense subset $U\subseteq Y$ where $\phi$ is a morphism. The
image, \cl J, of the induced homomorphism
$$\xymatrix{\O_Y^{N}\otimes\cl L\ve\ar@{->}[r]
\ar@{->>}[dr]&
  \O_Y
\\
&
\cl J\ar@{^(->}[u]\vphantom{I_{I_I}^I}
}
$$
is the sheaf of ideals defining the base locus. If $D$
denotes an effective Cartier divisor such that \,$\cl
J=\O_Y(-D)\cdot \cl J'$ \,for some ideal sheaf $\cl J'$,
then the set of zeros, $V(\cl J')$ is contained in $V(\cl
J)$. Clearly $\phi$ extends to the complement $U'=Y\setminus
V(\cl J')\supseteq U$ in such a way that the pullback of the
hyperplane bundle is
\\\centerline{$\phi_{|U}\us\O_{\ps{\C^N}}(1)=\cl
L\otimes\O(-D)$.}
\em\end{remark}

\subsection{(2,2,2)}

\medskip

When $q=2$ and $d_0=d_1=d_2=2$ the situation is still
manageable. It turns out that the indeterminacy locus of the
rational map
$$\ba{ccc}
\wt\rho: X=G(3,\s 2)&\rar&
\wt{\mathcal R}(r,\bar d)
\subset
\p{}(\s 3\otimes\wed3\s 1\us)\\ F_0\wedge
F_1\wedge F_{2} &\mapsto& dF_0\wedge dF_1\wedge dF_{2} \ea$$
is schematically equal to the image of the Veronese-like
embedding
$$
\ba{ccc}
Y=G(2,\s 1)&\stackrel \v\emb& X=G(3,\mathbf S_2) \\
\id{L_0,L_1}&\mapstochar\longrightarrow&\id{L_0^2, L_0L_1,
L_1^2}
.\ea
$$
Thus a single blowup $\pi:\wt X\ra X$ along $Y$
resolves the indeterminacy \ie the induced map $\wt\rho:\wt
X\ra\wt{\mathcal R}(r,\bar d)$ is a morphism. Indeed, write
the tautological sequence of $G(3,\s 2)$
\be\label{g3s2}
R_2\inj \s 2\surj Q_2
\ee
and likewise for $G(2,\s 1)$,
\be\label{g2s1}
R_1\inj \s 1\surj Q_1 .
\ee
The fiber of $R_2$ over $\u F\in X$ is
the space \id{F_0,F_1,F_{2}} spanned by three independent
quadratic forms.  In order to find the pullback of the
hyperplane class via the resolved map
$$
\xymatrix{
\wt X\ar@{->}[r]\ar@/^3.5mm/[rr]^{\wt\rho\pi} & X
\ar@{-->}[r]&
\wt{\mathcal R}(r,\bar d),}
$$
we have at first
$$
\xymatrix{
\wt\rho\us\O(-1)=\wed3R_2\ar@{->}[r]\ar@{->}[dr]
&\wed3\s 2 \ar@{->}[d] &\ni&
F_0\wed{}F_1\wed{}F_2 \ar@{|->}[d] \\
&\mathbf
S_3\otimes\wed3\s 1\us &\ni&
\,dF_0\wed{}dF_1\wed{}dF_2 .}
$$
The indeterminacy locus, $Z\subset X$, of $\wt\rho: X \rar
\wt{\mathcal R}(r,\bar d)$ is the scheme of zeros of the
slant arrow, $ \wed3R_2\lar \s 3\otimes\wed3\s1\us $.
Dualizing, we find $ \wed3 R_2\us \longleftarrow (\s 3
\otimes \wed3 \s 1\us)\us $, whence the ideal sheaf of $Z$
appears as the image
\be\label{Z}
\xymatrix{ (\s 3\otimes\wed3\mathbf
S_1\us)\us\otimes\wed3R_2 \ar@{->>}[r]& I(Z)\subset\O_X. }
\ee
We claim that $Z$ is equal to the image of $\v:G(2,\mathbf
S_1)\emb G(3,\s 2)$. Indeed, first note that $Z$ is
invariant under linear change of coordinates in \p r. Since
it is closed, it must contain a closed orbit of $G(3,\mathbf
S_2)$.  There are just two closed orbits, to wit those given
by the representatives: \id{x_0^2, x_0x_1, x_0x_2} and
\id{x_0^2, x_0x_1, x_1^2}. Only the latter one lies in
$Z$. The calculation of the tangent space to $Z$ at the
point \id{x_0^2,x_0x_1,x_1^2} performed below shows that $Z$
is of dimension at most $2(r-1)$. Since $Z$ contains the
image of $G(2,\s 1)$, it is in fact smooth and equal
to that image.  The tangent space is given by the equation
$$
\ba c
d(x_0^2+\varepsilon F_0)\wed{}d(x_0x_1+\varepsilon
F_1)\wed{}d(x_1^2+\varepsilon F_2)= \\
2\varepsilon dx_0\wed{}dx_1\wed{}( x_0^2d F_2-
2x_0x_1d F_1+ x_1^2d F_0 )=0,
\ea
$$
where the $F_i\in \s 2\big/\id{x_0^2,x_0x_1,x_1^2}$.

Equivalently:
$$
\left\{\ba c
\ds{ x_0^2\frac{\partial F_2}{\partial x_2}-
2x_0x_1\frac{\partial F_1}{\partial x_2}+
x_1^2\frac{\partial F_0}{\partial x_2}= \frac{\partial
}{\partial x_2}(x_0^2F_2- 2x_0x_1F_1+ x_1^2 F_0)=0,}
\\
\vdots\hskip1cm \\
\ds{ x_0^2\frac{\partial F_2}{\partial
x_r}- 2x_0x_1\frac{\partial F_1}{\partial x_r}+
x_1^2\frac{\partial F_0}{\partial x_r}= \frac{\partial
}{\partial x_r}(x_0^2F_2- 2x_0x_1F_1+ x_1^2 F_0)=0.}
\ea\right.
$$ We'd like to deduce that the subspace consisting of
triples
$$
(F_0,F_1,F_2)\in \left(\mathbf
S_2\big/\id{x_0^2,x_0x_1,x_1^2}\right)^{\oplus3}
$$
defined by the system just above must be of dimension
$$
\dim G(2,\s 1)=2(r-1).
$$
We see that $x_0^2F_2-2x_0x_1F_1+ x_1^2 F_0 $ is independent
of $x_2,\dots,x_r$. Thus, no monomial $x_mx_n,\,2\leq
m,n\leq r$ appears in the $F_i$. It follows that the $F_i$
are of the form
$$
F_i= a_{i0}x_0+a_{i1}x_1
$$ with the $a_{ij}\in\C[x_2,\dots,x_r]$ homogeneous of
degree one. We have then
$$\ba c
x_0^2(a_{20}x_0+a_{21}x_1) - 2x_0x_1(a_{10}x_0+a_{11}x_1) +
x_1^2 (a_{00}x_0+a_{01}x_1)= \\
a_{20}x_0^3+
(a_{21}-2a_{10})x_0^2x_1+
(a_{00}-2a_{11})x_0x_1^2 +a_{01}x_1^3
\in\C[x_0,x_1].
\ea$$ This implies
$$
a_{20}=a_{21}-2a_{10}=a_{00}-2a_{11}=a_{01}=0.
$$
Hence the $F_i$ depend exactly on $2(r-1)$
parameters. This achieves the verification that
$Z=\v(G(2,\s 1))$.
\medskip

Pulling back the surjection \rf{Z} to the blowup $\pi:\wt
X\ra X$, we find the surjections
$$
\xymatrix{\pi\us (\s 3\otimes\wed3\mathbf
S_1\us)\us\otimes\wed3R_2
\ar@{->>}[r]&\pi\us{}I(Z)\ar@{->>}[r] &\O_{\wt X}(-E)=I(E),
}
$$
with $E=\pi\inv Z$, the exceptional divisor. This yields
the formula
$$
\wt\rho\O_{\wt{\bbb F}}(1)= \pi\us\wed3R_2
\us\otimes\O_{\wt X}(-E).
$$
It follows that the pullback of the hyperplane class is
given by
$$
\wt\rho\us \h=\pi\us \q_1-E,
$$
where $\q_1=c_1Q_2$ (see\,\ref{g3s2}). Since \wt\rho\ is
generically injective, the degree of the image can be
calculated as
$$
\deg\mathscr R(r,2,2,2)=
\int_{\wt X}\wt\rho\us \h^{\dim X}.
$$
Setting $N=\dim X=\dim G(3,\s 2) =
3(\binom{r+2}{2}-3)$, we see that the degree is given
by
$$
\ba c\ds{ \int_{\wt X} \wt\rho\us \h^{N}= \int_X
\pi\ls\sum_0^{N} \binom{N}{i}\pi\us \q_1^i\cdot(-E)^{N-i}. }
\ea
$$
Using projection formula, we are reduced to the
calculation of
\bi
\item
the Pl\"ucker's degree of $G(3,\s 2)$ for the term
with $i=N$,
\ei and
\bi\item
the contribution of
{$\pi\ls(E)^j =(-1)^{j-1}\v\ls s_{j-\delta}\cl N, $}
 \ei
where $\cl N$ stands for the normal bundle of the embedding
\,$\v$\, and
$$
\delta=\rank \cl N=\dim G(3,\s 2)-\dim G(2,\mathbf
S_1).
$$
The minus signs come from the formula
$$
\iota\us\O_{\wt X}(E)=\O_{\cl N}(-1).
$$ The Segre classes of the normal bundle are obtained from
the usual exact sequence \be\label{N} \ba{*5c}
TY&\inj&TX_{|Y}&\surj&\cl N\\ ||& &||&&\\
\Hom(R_2,Q_2)&&\v\us\Hom(R_3,Q_3) && \ea\ee By definition of
\,$\v$, we have $\v\us R_3=\Sym_2R_2$.
Using \schub, we find,

\medskip

\centerline{
\begin{tabular}{c}
\begin{tabular}{|l|l|}
\hline $r$&{$\deg$}\\\hline 3&{1324220}\\\hline
4&{2860923458080}\\\hline 5&{243661972980477736263}\\\hline 6&
{728440733705107831789517245858}\\\hline 7&
{704613096513585123585398408696231899176183} \\
\hline
\end{tabular}\\
\framebox{$d_0=d_1=d_2=2$}
\end{tabular}
}

\medskip

\n A maple script is available at \cite{degsfol}.

\subsection{Bundles of projective spaces}

\medskip

When $k=2$ and $n_2=1$, the variety $X$ constructed in
\S \ref{S:X} is an open subset of a projective bundle over
an open subset of a grassmannian. In general we do not know
a manageable compactification.  Even when we can compactify
$X$ as above, the scheme structure of the base locus of
$\bar \rho$ can be non reduced and is far form being
understood in general.

Nevertheless in the following three cases  we are able to
handle the degree:
\bi\item
$q=1$ and $d_0$ divides $d_1$.
\item
arbitrary $q$ but
$k=2$ and $d_1=1$, \ie $\bar d=(1,\dots,1,e)$.
\item
 $q=1$, $d_0=2$ and $d_1=3$.
\ei
\subsubsection{
First Case: $q=1$ and  $d_0$ divides $d_1$.} This is in fact the
only case for which we got a closed formula. Now the natural
parameter space is the projective bundle
$$ X\lar\ps{\s {d_0}}
$$ described in the sequel.

Write the tautologic line subbbundle over $\ps{\mathbf
S_{d_0}}$,
$$ \O_{\s {d_0}}(-1) \inj \s {d_0}.
$$ Set $\kappa=d_1/d_0$. Taking symmetric power, we have the
exact sequence
$$ \O_{\s {d_0}}(-\kappa) \inj \mathbf
S_{d_1}\surj\G,
$$ which defines the vector bundle \G.  The fiber of \G\
over each $F_0 \in \ps{\s {d_0}}$ is the quotient
vector space $\s {d_1}/\id{F_0^\kappa}$.  Thus we
have
$$
\begin{array}{rcl}
\wt\rho:X=\ps{\G}&\lar& \wt{\mathscr R}(r, d_0, d_1)
\subseteq\ps{\s{d_1+d_0-2}\otimes \wed2\s 1\us}.
\\
(F_0,\ov F_1)&\mapstochar\longrightarrow& dF_0\wedge dF_1.
\end{array}
$$
The pullback of the hyperplane class via the map $\wt\rho$
is obtained as follows. Form the diagram
\be\label{spsq}
\xymatrix{
\O_{\s{d_0}}(-1)\otimes\s{d_1} \
\ar@{->}[dr]_\aa \ar@{>->}[r]& \s{d_0}\otimes \s{d_1}
\ar@{->}[d]&
\\ &\s{{d_1}+{d_0}-2}\otimes \wed2\s1\us
}
\ee
where the vertical map is defined by
$$
F_0\otimes F_1\mapsto dF_0\wed{}dF_1.
$$
Composing the slant arrow \,\aa\, with the natural
homomorphism
$$
\O_{\s{d_0}}(-1)\otimes \O_{\s {d_0}}(-\kappa)  \inj
 \O_{\s {d_0}}(-1)\otimes\s {d_1}
$$
we get zero since $dF_0\wed{}d(F_0^\kappa)=0$. Hence \,\aa\,
passes to the quotient,
$$ \xymatrix{
\O_{\s {d_0}}(-1)\otimes\s {d_1}
\ar@/^0.4cm/[rr]^\aa\ar@{->>}[r]&
\O_{\s {d_0}}(-1)\otimes\G{} \ar@{->}[r]_{\bar\aa}&
\s {{d_1}+{d_0}-2}\otimes \wed2\s 1\us.
}
$$
Composing \,$\bar\aa$\, with
$$
\O_{\s {d_0}}(-1)\otimes \O_{\G}(-1)  \inj
\O_{\s {d_0}}(-1)\otimes\G
$$
we finally find the line subbundle,
$$
\O_{\s {d_0}}(-1) \otimes\O_{\G}(-1)\inj \s{{d_1}+{d_0}-2}
\otimes\wed2\s 1\us.
$$
The last map is injective at the point
$(x_1^{d_0},\ov{x_1^{d_1-1}x_2})$, which is a representative
of the unique closed orbit of \ps{\G}. Hence it is injective
everywhere.  Alternatively, since $F_0^\kappa,F_1$ are
linearly independent, the rational map $\p r\rar\p1$ they
define is non-constant, hence $dF_0\wed{}dF_1\nn0$. Thus,
the pullback to $X$ of the hyperplane class of the
projective space
$\p{}\big(\s {{d_1}+{d_0}-2}\otimes\wed2\s1\us\big)$
is
$$
H=\h+\h'
$$
where $\h=c_1\O_{\s {d_0}}(1)$, which comes from the
base $\ps{\s {d_0}}$, \ and \
$\h'=c_1\O_{\G}(1)$, the relative hyperplane class. With the
notation as in \rf{q1..1}, we have
$$
\rank\G-1=N_{d_1} -2
$$
for the fiber dimension of $\ps\G\ra\ps{\mathbf
S_{d_0}}$. The sought for degree is
\be\label{pq}\ba{rl}
\deg \mathscr R(r, d_0,d_1)=&
\ds{\int_{\ps\G}}H^{N_{d_1}+N_{d_0}-1} =\ds{\sum_{i}}
\bsm{N_{d_1} +N_{d_0} -1}\\ {i}\esm \h^i s_{N_{d_0}-i}(\G)\\\na3
=&\ds{\binom{N_{d_1} +N_{d_0} -1}{N_{d_0} }-
\frac{d_1}{d_0}\binom{N_{d_1} +N_{d_0} -1}{N_{d_0} -1}}
.\ea
\ee
The last equality follows from the calculation of the Segre
class $s(\G)=1-\kappa \h$, so $s_i(\G)$ is zero in degrees
$i\geq2$.

If $r=3,\,{d_1}=2,\,{d_0} =1$, one finds
$\binom{3+8}3-2\binom{11}2=55$. By constrast, the degree of
the Segre variety $\pd3\vez\p9\subset\p{39}$ of which the
image of $\rho$ is a rational projection, is equal to
$\binom{12}{3}$.

\subsubsection{Second case: $k=2$ and $d_0=1$.}

We are now looking at foliations defined by
$\omega=i_R(dF_0\wed{}\cdots\wed{}dF_q)$ where $\deg
F_0=\cdots=\deg F_{q-1}=1;\,\deg F_q=d\geq2$. A natural
parameter space is the projective bundle over the
grassmannian $G=G(q,\s 1)$ defined as follows. Write
the tautological sequence
$$ R_q\inj \s 1\surj Q.
$$ The fiber of $R_q$ over $\u F\in G$ is the space
\id{F_0,\dots,F_{q-1}} spanned by linear forms. Now the last
polynomial $F_q$ is taken as a class in the projective space
\ps{\s {d}/\id{F_0^d,F_0\cdot{}F_1^{d-1},\dots,
F_{q-1}^d}}.  The natural homomorphism $\Sym_dR_q\ra \mathbf
S_d$ is injective; it corresponds to an instance of the
vector bundle $\cl A_2$ described  in \ref{S:X}. Form the
projective bundle
$$
\pi:X=\ps{\s {d}/\Sym_dR_q}\lar G.
$$
 Note that the rational map
$$
\ba{ccc} X&\xymatrix{\ar@{-->}[rr]^{\bar\rho}&&}&
\p{}(\s{d-1}\otimes\wed q\s 1\us) \\
(\id{F_0,\dots,F_{q-1}},\ov F_q)
&\xymatrix{\ar@{|->}[rr]&&}&
dF_0\wedge\dots\wedge{}dF_{q-1}\wedge{}dF_q \ea$$ is in fact
regular everywhere. Indeed, regularity is an open condition;
the map is invariant under the natural action of $GL_{r+1}$
and is regular at the representative
$(\id{x_0,\dots,x_{q-1}},\ov {x_q^{d-1}x_0})$ of the unique
closed orbit. Thus the sought for degree can be computed by
Schubert calculus in the following manner. Set
\be\label{Ng1}
\ba l g=q(r+1-q)=\dim G\\
N=\mbox{$\binom{r+d}{r}$}-\mbox{$\binom{q-1+d}{q-1}$}-1,
\ea
\ee
so that presently the dimension of the component is
$\delta=N+g$. The pullback of the hyperplane class from
$\p{}(\s {d-1}\otimes\wed q\s 1\us)$ is
equal to $\h+\q_1$,
where $\h$ stands for the relative hyperplane class of the
projective bundle $X\ra G$\, and\, $\q_1=c_1Q$. By general
principles, the degree is given by
$$
\int_X(\h+\q_1)^\delta=\sum_0^g\binom{\delta}{i} \int_G
\pi\ls(\h^{\delta-i})\q_1^i =\sum_0^g\binom{\delta}{i}
\int_G s_{g-i}\cdot\q_1^i.
$$
Here $s_i=c_i(\Sym_dR)$. For $q=2,\,r=3$ we find
$$
d^2(d-1)(d+3)(d^2+2)(d^2+4d+6)(d+2)^2(d+1)^2
\big/(2^6\cdot3^5),
$$
a polynomial of degree 12 in $d$. For $q=2;\,r=4,5,6,7,8$ we find
polynomial formulas of respective degrees 24,\,40,\,60,\,84,\,112.
This suggests a polynomial degree like $2r(r-1)$. Now for
$q=3,\,r=4,5,6,7,8$ we get polynomial formulae of degrees $3r(r-2)$
with respect to $d$. Further experiments (cf.\,\cite{degsfol})
suggest polynomial formulas of degrees $qr(r-q+1)$. Here is a sample
for small values of $r,q,d$.
$$\ba c \framebox{$(r,q)=(5,2)$}\\
\ba{|c|c|c|c|c|}\hline
d&2&3&4&5\\\hline
\deg&2390850&
10457430102&
9654013512864&
3099059696318355\\
\hline
\ea
\ea$$
{\footnotesize
$$\ba c \framebox{$(r,q)=(6,2)$}\\
\ba{|c|c|c|c|c|}\hline
d&2&3&4&5\\\hline
\deg&1139133688
&91451421683006
& 1118409272891730904
& 3524857658574891999976

\\
\hline
\ea
\ea$$
$$\ba c \framebox{$(r,q)=(6,3)$}\\
\ba{|c|c|c|c|}\hline
2&3&4&5\\\hline
 8983484048
&9350781792221835
&1060759743612735149417
&22044166363067583367287424
\\
\hline
\ea
\ea$$
}

\subsection{${(2,2m+1)}$}

Assume $q=1$, $d_0=2$ and $d_1=3$. \label{SS:23} Set for short
$\mathbf X=\ps{\s2}\vez\ps{\s3}$. Put as before $N_d =
\binom{r+d}d-1$. We have
$$
 \dim \x =N_2+N_3.
$$
We look closer at the indeterminacy locus of
$$\ba{ccc}
\wt\rho:\x  &\dasharrow &
\ps{\s3\otimes\wed2\s1\us} \\
(F,G) & \mapsto & dF\wed{}dG.
\ea$$
It is, set-theoretically,
$$
\B(\wt\rho)=\{(L^2,L^3)\st
L\in\ps{\s1}\}.
$$

\begin{lemma}
The tangent space to the scheme of
indeterminacy $\mathbf{B}=\mathbf{B}(\wt\rho)$
is the subspace
$$
\{ (F',G')\in T_{(L^2,L^3)}\mathbf X=\s2/\id{L^2}\oplus\s3/\id{L^3} \st G'=\frac32 F' \}.
$$
\end{lemma}

\begin{proof}
The tangent space to the scheme of
indeterminacy is the set of pairs $(F',G')$
such that $d(L^2+\varepsilon F')\wedge
d(L^3+\varepsilon G')=0 $.  Expanding we get
\be\label{tgb}
2dL\wedge dG'+ 3LdF'\wedge dL
=dL\wedge (2dG'- 3LdF') =0 .
\ee
By division, we must have $2dG'- 3LdF'=F''dL$
for some $F''\in\s2$.  This implies $dF''\wedge
dL=3dF'\wedge dL$. Hence again by division,
$dF''-3dF'=A dL$ for some $A\in \s1$. This
implies $A=aL$ for some constant $c$. Thus
$d(F''-3F'-\frac12aL^2)=0$ so that in fact
$F''=3F'+\frac12aL^2$. Plugging back in a
previous relation, we find $2dG'-
3LdF'=(3F'+\frac12aL^2)dL$ whence
$2dG'-\frac16adL^3= 3d(LF')$. This yields $2G'
-\frac16aL^3= 3LF'$, hence $G'=\frac32LF'$ in
$\s3/\id{L^3}$. Conversely, it is easy to see
that for such $G'=\frac32LF'$, the differential
form $2dG'- 3LdF'$ is a
multiple of $dL$, hence \rf{tgb} holds.
\end{proof}

\medskip
Set $\V =\B_{red} \cong
\ps{\s1}$. Thus $\B$ is a multiple
structure or thickenning of $\V $. The
tangent sheaf to $\B$ is in fact a
vector bundle of rank $\dim\ps{\s2}$. We have
the exact sequence of vector bundles over
$\V $,
$$
T\V \inj T\B_{|\V }\surj
\cl N_{\V /\B}
$$
where $\cl N_{\V /\B}$ stands
for the normal bundle of $\V \subset
\B$. We register the formula
$$
\rank\cl N_{\V /\B}
=\binom{r+2}{2}-r=\binom{r+1}{2}+1.
$$
We look at the blowup $\x '\ra \mathbf
X$ along $\V $. Denote by $\mathbf
E\subset \x '$ the exceptional divisor.
Recall we have $\EE'=\ps{\cl N_{\mathbf
V/\x }}$, the projectivization of the
normal bundle of $\V \subset \x $

\begin{lemma}\label{ro1} We assume $r \le 5$.
Let $\rho': \x '\rar \ps{\s3\otimes\wed2\s1\us}$ be the rational map
induced by $\wt\rho$ and denote by  $\B' \subset \x '$
the indeterminacy scheme of $\rho'$.  Then we have
$$\ba c \B'=\ps{\cl N_{\V /\mathbf B}}\subset \ps{\cl N_{\V /\x}}=\EE', \ea$$
the projectivization of the normal bundle of $\V $ in its thickenning $\B$.
\end{lemma}

\begin{proof}
We look at the diagram of tangent/normal
bundles over $\V $,
\be\label{nvb}\xymatrix{ \ T\V  \
\ar@{=}[r]\vphantom{I_{I_I}^I}\ar@{>->}[d]& \
T\V  \
\vphantom{I_{I_I}^I}\ar@{>->}[d]&\\ T\mathbf
B_{|\V } \ar@{->>}[d] \
\ar@{>->}[r]&T\x _ {|\V }
\ar@{->>}[d]\ar@{->>}[r] &\cl N_{\mathbf
B/\x |\V }\ar@{=}[d] \\ \cl
N_{\V /\B} \ \ar@{>->}[r]&\cl
N_{\V /\x } \ar@{->>}[r]& \cl
N_{\B/\x |\V } }
\ee
which tells us that \ps{\cl N_{\V/\B}} embeds naturally into
$\EE' = \ps{\cl N_{\V /\x }}$. Let $x'\in\EE'$. Thus we may
represent it as $x'=\lim_{\varepsilon \ra0} (L^2 +
\varepsilon F',L^3+\varepsilon G')$ for some $(F',G') \in
T_{(L^2,L^3)} \x $ with nonzero image in $\cl N_{\V / \x
}$. Here we think of $(L^2 + \varepsilon F', L^3+
\varepsilon G')$ as a small arc in $\x \setminus \V $ for
$\varepsilon\nn0$. Hence it lifts to an arc in $ \x '
\setminus \EE'$ which hits $x'\in \EE'$ for
$\varepsilon=0$. As in \rf{tgb} we find for
$\varepsilon\nn0$,
\be\label{df1dg1}\ba{rl}
\rho(L^2 + \varepsilon F', L^3+
\varepsilon G') & =\varepsilon LdL\wedge (2dG'-
3LdF') + \varepsilon^2dF' \wedge dG'\\
&=LdL\wedge (2dG'- 3LdF') + \varepsilon dF'
\wedge dG'.
\ea\ee
Now if $x'$ is {\em not} in the indeterminacy locus, $\B'$,
then we must have
$$
\rho'(x')=\lim_{\varepsilon \ra0}\rho(L^2 +
\varepsilon F', L^3+ \varepsilon G').
$$
This limit is $\rho'(x')=LdL\wedge (2dG'-
3LdF')$ provided the expression is \nn0. It is zero if and
only if  $G'=\frac23LF'$, \ie  $x'$ is  in
$\ps{\cl N_{\V /\B}}$. In this case, recalling\,\rf{df1dg1},
$$ \rho'(x')=dF' \wedge
dG'=\frac23F'dF'\wedge{}dL.
$$
Since the right hand side must be
(projectively) independent of  representatives
of $F'\in\s2/\id{L^2}$, we must have $dL\wedge
dF'=0$, a contradiction. Thus $LdL\wedge (2dG'-
3LdF') $ must be \nn0, \ie $x'$ is {\em not} in
$\ps{\cl N_{\V /\B}}$.  This yields
$ \ps{\cl N_{\V /\B}}
\subseteq\B' $.

The dimension is given by
$$\ba c
\dim\B'= \dim\V +\rank\cl
N_{\V / \B}-1=
r+\binom{r+2}2-1-r-1 =N_2-1=\binom{r+2}2-2.
\ea$$
Thus we also have $\cod\B'=\rank \cl N_{\B'/\x'}=N_3+1$.

Unfortunately, for the other inclusion
we don't know how to proceed coordinate-freewise.
Using coordinates, with the help of computer algebra ({\sc singular}),
it can be checked (see \cite{degsfol})  that $\B'$ is smooth and of the right
dimension $\dim\ps{\cl N_{\V /\B}}$.  This requires fixing $r$ to low values,
\eg $r\leq5$.  Here is an outline of the calculation for $r=2$.
We take affine coordinates
$a_1,\dots,a_5,b_1,\dots,b_9$ for \ps{\s2}\vez\ps{\s3}. Set
$$
\ba l F=x_{0}^{2} +a_{1}x_{0}x_{1}
+a_{2}x_{0}x_{2}+a_{3}x_{1}^{2}
+a_{4}x_{1}x_{2} +a_{5}x_{2}^{2}, \\
G=x_{0}^{3} +b_{1}x_{0}^{2}x_{1} +b_{2}x_{0}^{2}x_{2}+\cdots
+b_{8}x_{1}x_{2}^{2} +b_{9}x_{2}^{3}.  \ea
$$
We compute $ dF\wedge{}dG $ expanding the 2\vez2 minors of
the 2\vez3 matrix with rows the gradients of $F,G$. We find
three cubics as coefficients of
$dx_0\wedge{}dx_1,dx_0\wedge{}dx_2,dx_1\wedge{}dx_2$. The
indeterminacy locus, \B, is given by the ideal spanned by those
thirty coefficients.  Its jet of order one is spanned by
nine independent linear equations, in agreement with the
expected tangent space dimension, to wit, 5, the freedom of
the quadric $F$. Continuing, we find next the local equations of the
bi-Veronese, eliminating $c_1,c_2$ from the 5+9
equations obtained from the conditions
$$
F=(x_0+c_1x_1+c_2x_2)^2,\,G=(x_0+c_1x_1+c_2x_2)^3.
$$
We find that the ideal of the bi-Veronese is spanned by the
12 polynomials
$$\ba c 2b_{1} -3a_{1},\, 4b_{3} -3a_{1}^{2},\,
8b_{6} -a_{1}^{3},\, 2b_{2} -3a_{2},\, 2b_{4}
-3a_{1}a_{2},\, 8b_{7} -3a_{1}^{2}a_{2}, \\
4b_{5} -3a_{2}^{2},\, 8b_{8}
-3a_{1}a_{2}^{2},\, 8b_{9} -a_{2}^{3},\,
a_{1}^{2} -4a_{3},\, a_{1}a_{2} -2a_{4},\,
a_{2}^{2} -4a_{5}.  \ea
$$
Accordingly, the blowup is covered by 12 affine patches, one
for each choice of the principal generator for the
exceptional ideal.  The 9 generators involving a
$b-$coefficient belong to the ideal of \B.  It follows that
the indeterminacy locus $\B'$ is disjoint from these nine
neighborhoods.  We are left with the 3 equations
$4a_{3}-a_{1}^{2} , 2a_{4}-a_{1}a_{2},4a_{5} -a_{2}^{2}$;
these define the Veronese in \p5.  Choosing
$\varepsilon=a_{1}^{2} -4a_{3}$ as the exceptional
generator, the blowup is written as
$$
\left\{\ba l

b_1= \frac{1}{2}\varepsilon c_{1} +\frac{3}{2}a_{1},\,
b_2= \frac{1}{2}\varepsilon c_{2} +\frac{3}{2}a_{2},\\\na3
b_3= \frac{1}{4}\varepsilon c_{3} +\frac{3}{4}a_{1}^{2},\,
b_4= \frac{1}{2}\varepsilon c_{4} +\frac{3}{2}a_{1}a_{2},\\\na3
b_5= \frac{1}{4}\varepsilon c_{5} +\frac{3}{4}a_{2}^{2},\,
b_6= \frac{1}{8}\varepsilon c_{6} +\frac{1}{8}a_{1}^{3},\\\na3
b_7= \frac{1}{8}\varepsilon c_{7}
+\frac{3}{8}a_{1}^{2}a_{2},\,
b_8= \frac{1}{8}\varepsilon c_{8}
+\frac{3}{8}a_{1}a_{2}^{2},\\\na3
b_9= \frac{1}{8}\varepsilon c_{9} +\frac{1}{8}a_{2}^{3},\,
a_4= -\frac{1}{2}\varepsilon c_{10}
+\frac{1}{2}a_{1}a_{2},\\\na3
a_5= -\frac{1}{4}\varepsilon c_{11} +\frac{1}{4}a_{2}^{2}.

\ea\right.
$$
Substituting into the ideal of the
indeterminacy locus, the original 30 generators
become divisible by the local equation,
$\varepsilon$, of the exceptional
ideal. Dividing, we obtain the ideal of the
indeterminacy locus upstairs,  that is, of the
induced rational map $\rho'$
(cf.\,Lemma\,\ref{ro1}). We find the ideal of $\B'$ is
presently generated by
$$\begin{array}{c}
 c_{1},c_{2}, 2c_{4} +3c_{10}, 2c_{7} +6c_{10}a_{1} -3a_{2},
2c_{5} +3c_{11}, \\
2c_{8} +6c_{10}a_{2} +3c_{11}a_{1},
2c_{9}
+3c_{11}a_{2},  2c_{3} -3, 2c_{6} -3a_{1},
\underbrace{a_{1}^{2}-4a_{3}}_\varepsilon.
\end{array}
$$
Thus we see that the indeterminacy locus is contained in the
exceptional divisor and we also learn that it is in fact a
projective subbundle of the exceptional divisor $\EE'$, in agreement with \ref{ro1}.
\end{proof}

\begin{remark}\rm
Lemma \ref{ro1} is valid only for small values of $r$,
as stated and explained in the proof.
But we conjecture that it is true for all $r$.
It seems that a more conceptual proof is needed and
probably  it would involve some new idea.
\end{remark}

At any rate,  for each value of $r$, the validity of  Lemma  \ref{ro1}
is all we need to find the degree:
We consider the following diagram displaying the
resolution of the map $\wt\rho: X \rar\ps{\s3\otimes\wed2\s1\us}$.
$$
\xymatrix@R-10pt{
\EE''\ar@{->}[d]
& & &\subset&\x ''
\ar@{->}[d]
\ar@{->}[ddr]^{\rho''}
_{\wt\rho'\hskip.197cm}
&
\\
\B'&\subset&
\EE'\ar@{->}[d] &\subset&\x '
\ar@{->}[d]
\ar@{-->}[dr]
&
\\&&
\V  &\subset&\x \ar@{-->}[r]^
{\wt\rho\hskip.97cm}&
\ps{\s3\otimes\wed2\s1\us}
}
$$
The pullback of the hyperplane class via
$\rho''$ can be written as
$$\rho''\inv\h=m_1\h_1+m_2\h_2+m_3\mathbf
e'+m_4\mathbf e''
$$
for suitable integers $m_i$, where we've denote the cycles
$\E'=[\EE'],\E''=[\EE'']$ and $\h_i$ the hyperplane class of
each factor in $\x=\ps{\s2}\vez\ps{\s3}$.  The coefficients
$m_i$ will be determined using the Remark\,\ref{blocus} and
excision (cf.\,\cite[1.8,\,p.\,21]{Fulton}). Over $\mathbf
U=\x \setminus\V $ only $\h_1,\h_2$ survive and we have
$\rho''\inv_{\mathbf U}\h= \wt\rho\inv_{\mathbf
U}\h=\h_1+\h_2 $ since $\rho_{\mathbf U}$ is defined by a
bihomogeneous expression of bidegree 1,1.  Put $\mathbf
U'=\x '\setminus \B' =\x ''\setminus \EE''$.  The local
calculations show that the image of
$(\s3\otimes\wed2\s1\us)\us \otimes \O_{\p{}(\s2)}(-1)
\otimes \O_{\p{}(\s3)}(-1)\ra\O_{\x'}$ is equal to
$\O_{\x'}(-\EE')\cdot\cl I(\B')$. Blowing-up $\B'$, we find
the surjection
$$
(\s3\otimes\wed2\s1\us)\us \otimes \O_{\p{}(\s2)}(-1)
\otimes \O_{\p{}(\s3)}(-1)\otimes
\O_{\x''}(\EE')\twoheadrightarrow
\O_{\x''}(-\EE'')=
\cl I(\B')\O_{\x''}
.
$$
Thus, we have
$$
(\wt\rho)\inv\h=\h_1+\h_2-\E'-\E''.
$$
The degree is computed as
$$
\int_{\x ''}(\h_1+\h_2-\E'-\E'')^{N_2+N_3}.
$$
Apart from the term
$\int_{\x ''}(\h_1+\h_2)^{N_2+N_3}=\binom{N_2+N_3}{N_2}
$, all others lie over $\V$. Since
$\h_1\cap\V=2\h,\h_2\cap\V=3\h$ and $\h^{r+1}=0$, we see
that terms like
$\h_1^i\h_1^j(\E')^k(\E'')^l$  give zero whenever $i+j>r$.
Thus the relevant part of the integrand is
$$
\sum_0^r\binom{N_2+N_3}{i}(5\h)^i(-\E'-\E'')^{N_2+N_3-i}.
$$
First we collect coefficients of $\E''$, then take the
pushforward to $\x '$ using our knowledge of the
normal bundle of $\B'\subset \x'$ and so on
till $\x $. Thus
$$
(\E'')^i=(\E'')^{i-1}\E''\leadsto
(-1)^{i-1}s_{i'}(\cl N_{\B'/\x'})\cap[\B'],
$$
with $i'=i-\cod\B'=i-N_3-1$.
The sum above pushes forward to
$$\ba{rl}
\ds{\sum_0^r}\mbox{$\binom{N_2+N_3}{i}$}
(5\h)^i(-1)^{N_2+N_3-i}
&\left(-(\E')^{N_2+N_3-i} \,+\right.\\
&
\ds{\sum_{j=\dim\B'}^{N_2+N_3-i}}\mbox{$
\binom{N_2+N_3-i}{j}$}
(\E')^{N_2+N_3-i-j}(-1)^{j-1}s_{j'}
\mbox{$\left.\vphantom{(I')_I^I}\right)$},
\ea$$
setting for short $s_{j'}=s_{j'}(\cl N_{\B'/\x'})\cap[\B']$,
with
$j'=j-\cod\B'=j-N_3-1$.
(Thus $s_{j'}$ is a cycle of dimension
$N_2-1-j'=N_2+N_3-j$.)
These Segre classes can be derived from \ref{ro1} as
follows. We have the exact sequence
\be\label{nb1x1}
\cl N_{\B'/\EE'} \inj \cl N_{\B'/\x'} \surj
\O_{\EE'}(\EE')_{|\B'}.
\ee
We also recall that, for any exact sequence of vector
bundles
$$
\cl E'\inj\cl E\surj\cl E''
$$
we have the formula for the normal bundle of
$\ps{\cl E'}\subset\ps{\cl E}$
$$
\cl N_{\ps{\cl E'}/\ps{\cl E}}=\cl E''\otimes\O_{\cl E'}(1).
$$
In view of \rf{nvb}, this yields
$$
\cl N_{\B'/\EE'}=
\cl N_{\B/\x|\V}\otimes \O_{\cl N_{\V/\B}}(1).
$$
The actual calculation is best performed using computer
algebra. A script using {\sing} is available at \cite{degsfol}.
A sample of the first few values is listed below.
$$
\ba c
\ba{|c|l|}
\hline
r&\deg\\\hline
2&770\\\hline

3&6254612\\\hline

4&481152797320\\\hline

5&803161672838504856\\\hline

6&36968358460592709286459400\\\hline

7&53639021695280557844870264612516640\\\hline

8&2759237622445467221610266591396121818496881016\\\hline

\ea
\\\na{.158}
\framebox{(2,3)}
\ea
$$

As a final remark we mention that there is compelling
computer algebra evidence indicating that the case of
bidegree (2,3) carries over to the case $(2,2m+1)$ with
slight modifications. The indeterminacy locus of the
rational map $\x=\ps{\s2}\vez\ps{\s{2m+1}} \rar
\p{}(\s{2m+1}\otimes\wed2\s1\us)$ given by $(F,G)\mapsto{}
dF\wedge{} dG$ is again a thickening of the biveronese
$\{(L^2,L^{2m+1})\st L\in\ps{\s1}\}$. Blowing up the reduced
structure, the indeterminacy locus, $\B'$, of the induced
rational map $\x'\rar \p{}(\s{2m+1}\otimes\wed2\s1\us)$ is
no longer reduced for $m>1$.  Nevertheless, it still is a
rather manageable complete intersection. In fact, we find
local equations of $\B'$ of the form $e^m,f_1,...,f_u$, with
$e$ denoting the equation of the exceptional divisor, and
the $f_i$'s define a projective subbundle of the exceptional
divisor just as in the case (2,3).

$$
\ba{c}
\framebox{$\deg\mathscr R(r,d_0=2,d_1=2m+1)$}\\
\ba{|c|l|}
\hline
d_1&\deg\ (\p3) \\\hline
5&27500627268\\\hline
7&19062120397608\\\hline
9&3910289698588916\\\hline
11&341013122932980120\\\hline

\ea
\qquad
\ba{|c|l|}
\hline
d_1&\deg \ (\p4)
\\\hline

5&5858652068789831804\\\hline
7&2734930355086609774678630\\\hline
9&118796991387599661786404269060\\\hline
11&955667356931740162987705236374200\\\hline

\ea

\ea
$$

Interpolating the first few values of odd $d_1$, we find for
\p3 the polynomial
\\
$( t-1 )  \big( t^{26}+55
t^{25}+1450t^{24}+24616t^{23}+305020t^{22}+2961172t^{
21}+23561656t^{20}+158392960t^{19}+918866662t^{18}+
4670514826t^{17}+21033417148t^{16}+84615935632t^{15}+
305921226844t^{14}+998318576836t^{13}+2949392111320t^{12}+
7903552056256t^{11}+19229223618721t^{10}+41774679574903t^{
9}+72390849730794t^{8}+15945324910344t^{7}-541088235621216{t
}^{6}-2539188961011216t^{5}-315410776482528t^{4}+
14933666207688192t^{3}+85822791395378688t^{2}-
247712474710388736t+162893498195312640
\big)/3656994324480
$.
\\It fits all values of $\deg\mathscr R(3,2,t),\,t=2m+1$,
up to $m=35,d_1=71$, presently the physical limit of our
computer's memory. It should be noted that $\deg\mathscr R(3,2,2t)=\binom{N_{2t} +N_{2} -1}{N_{2} }-
\frac{2t}{2}\binom{N_{2t} +N_{2} -1}{N_{2} -1}$
is  a polynomial in $t$ of the same degree 27 as above.

\vspace{3 cm}

\begin{tabular}{lll}
 {\small Fernando Cukierman}& \hbox {\small Jorge Vit\'orio Pereira }& \hbox {\small Israel Vainsencher} \\
 {\small Depto. Matem\'atica, FCEN-UBA} & \hbox{\small IMPA} & \hbox {\small  Depto. Matem\'atica, UFMG} \\
 {\small Ciudad Universitaria} & \hbox{\small Estrada Dona Castorina 110} & \hbox{\small Av. Antonio Carlos 6627} \\
 {\small 1428 Buenos Aires}& \hbox{\small 22 460-320 Rio de Janeiro}& \hbox{\small 31 270-901 Belo Horizonte} \\
 {\small Argentina}& \hbox{\small Brasil}& \hbox{\small Brasil} \\
 & & \\
 {\small fcukier@dm.uba.ar}& {\small jvp@impa.br}& {\small israel@mat.ufmg.br}
 \end{tabular}

\end{document}